\documentstyle[12pt]{article}

\newtheorem{definition}{Definition}[section]
\newtheorem{theorem}[definition]{Theorem}
\newtheorem{lemma}[definition]{Lemma}
\newtheorem{corollary}[definition]{Corollary}
\newtheorem{remark}[definition]{Remark}
\newtheorem{example}[definition]{Example}
\newtheorem{conjecture}[definition]{Conjecture}
\newtheorem{problem}[definition]{Problem}
\newtheorem{note}[definition]{Note}

\def\C{{\rm I}\!\!\! {\rm C}}

\newcommand{\fld}{{\cal F}}

\newcommand{\ls}{{\Phi = (A;
 E_0,E_1,\ldots,E_d;
A^*; E^*_0,E^*_1,\ldots,E^*_d)}}
\newcommand{\Span}[1]{{ \mbox{\rm Span}\{#1\}}} 
\newcommand{\beast}{\begin{eqnarray*}}
\newcommand{\eeast}{\end{eqnarray*}}

\begin{document}
\newenvironment{proof}{\noindent{\it Proof\/}:}{\par\noindent $\Box$\par}

\title{ \bf Some algebra related to $P$-and $Q$-polynomial \\
association schemes}
\author{Tatsuro Ito,  
Kenichiro Tanabe, and  
Paul Terwilliger
}
\date{}
\maketitle
\begin{abstract}

Inspired by the theory of $P$-and $Q$-polynomial association
schemes we consider the following situation in linear algebra.
Let $\fld$ denote a field, and let $V$ denote
a vector space over $\fld$ with finite positive dimension.
We consider a pair  of linear transformations  
  $A:V\rightarrow V$ and  $A^*:V\rightarrow V$
  satisfying
the following four conditions.
\begin{enumerate}
\item $A$ and $A^*$ are both diagonalizable on $V$.
\item There exists an ordering $V_0, V_1,\ldots, V_d$ of the  
eigenspaces of $A$ such that 
\beast
A^* V_i \subseteq V_{i-1} + V_i+ V_{i+1} \qquad \qquad (0 \leq i \leq d),
\eeast
where $V_{-1} = 0$, $V_{d+1}= 0$.
\item There exists an ordering $V^*_0, V^*_1,\ldots, V^*_\delta$ of the 
eigenspaces of $A^*$ such that 
\beast
A V^*_i \subseteq V^*_{i-1} + V^*_i+ V^*_{i+1} \qquad \qquad (0 \leq i \leq \delta),
\eeast
where $V^*_{-1} = 0$, $V^*_{\delta+1}= 0$.
\item There is no subspace $W$ of $V$ such  that  both $AW\subseteq W$,
$A^*W\subseteq W$, other than $W=0$ and $W=V$.
\end{enumerate}
We call such a pair a TD pair.
Referring to the above TD pair, 
we show $d=
\delta$. We show that for $0 \leq i \leq d$, the eigenspaces
$V_i$ and $V^*_i$ have the same dimension. Denoting this common dimension by
$\rho_i$, we show the sequence $\rho_0,\rho_1,\ldots, \rho_d$ is symmetric and
unimodal, i.e. $\rho_{i-1} \leq \rho_i$ for $1 \leq i \leq d/2$ and 
$\rho_{i} = \rho_{d-i}$ for $0 \leq i \leq d$. We show that there exists 
a sequence of scalars $\beta, \gamma, \gamma^*, \varrho, \varrho^*$ taken from
$\fld$ such that both 
\beast
0 &=&\lbrack A,A^2A^*-\beta AA^*A + 
A^*A^2 -\gamma (AA^*+A^*A)-\varrho A^*\rbrack,
\\
0 &=& \lbrack A^*,A^{*2}A-\beta A^*AA^* + AA^{*2} -\gamma^* (A^*A+AA^*)-
\varrho^* A\rbrack,  
\eeast
where $\lbrack r,s\rbrack $ means $rs-sr$. 
The sequence is unique if $d\geq 3$. 
Let $\theta_i$ (resp. $\theta^*_i$) denote the eigenvalue of $A$ (resp. $A^*$)
associated with $V_i$  (resp. $V^*_i$), for $0 \leq i \leq d$. We show the
expressions 
\beast
{{\theta_{i-2}-\theta_{i+1}}\over {\theta_{i-1}-\theta_i}},\qquad \qquad  
 {{\theta^*_{i-2}-\theta^*_{i+1}}\over {\theta^*_{i-1}-\theta^*_i}} 
 \qquad  \qquad 
\eeast 
 both equal $\beta +1$, for $\;2\leq i \leq d-1$.  
We hope these results will ultimately lead to a complete classification
of the TD pairs.

\end{abstract}


%
%
%
%
%
%
%
%


\section{Introduction}


\medskip
\noindent 
We begin with the following situation in linear algebra.
Let $\fld$ denote a field, and let $V$ denote
a vector space over $\fld$ with finite positive dimension.
By a {\it linear transformation on} $V$,
we mean a $\fld$-linear
map from $V$ to $V$. Let $A$ denote a linear transformation
on $V$. By an {\it eigenspace} of $A$, we mean 
a nonzero  subspace of $V$ of the form
\begin{equation}
\lbrace v \in V \;\vert \;Av = \theta v\rbrace,
\label{eq:defeigspace}
\end{equation}
where $\theta \in \fld$.
We say $A$ is {\it diagonalizable} on $V$ whenever
$V$ is spanned by the eigenspaces of $A$.

\begin{definition}
\label{def:nonthinlp}
Let $\fld$ denote a field, and let $V$ denote
a vector space over $\fld$ with finite positive dimension.
By a Tridiagonal pair (or TD pair) on $V$,
we mean an ordered pair $(A, A^*)$, where
$A$ and $A^*$ are linear transformations on $V$ that satisfy
the following four conditions.
\begin{enumerate}
\item $A$ and $A^*$ are both diagonalizable on $V$.
\item There exists an ordering $V_0, V_1,\ldots, V_d$ of the  
eigenspaces of $A$ such that 
\begin{equation}
A^* V_i \subseteq V_{i-1} + V_i+ V_{i+1} \qquad \qquad (0 \leq i \leq d),
\label{eq:astaractionthreeint}
\end{equation}
where $V_{-1} = 0$, $V_{d+1}= 0$.
\item There exists an ordering $V^*_0, V^*_1,\ldots, V^*_\delta$ of
the  
eigenspaces of $A^*$ such that 
\begin{equation}
A V^*_i \subseteq V^*_{i-1} + V^*_i+ V^*_{i+1} \qquad \qquad (0 \leq i \leq \delta),
\label{eq:aactionthreeint}
\end{equation}
where $V^*_{-1} = 0$, $V^*_{\delta+1}= 0$.
\item There is no subspace $W$ of $V$ such  that  both $AW\subseteq W$,
$A^*W\subseteq W$, other than $W=0$ and $W=V$.
\end{enumerate}
\end{definition}
\begin{note} According to a common 
notational convention, for a linear transformation $A$
the 
conjugate-transpose of $A$ is denoted
$A^*$.
We emphasize we are {\it not} using this convention.
In a TD pair $(A,A^*)$, the linear transformations 
$A$ and $A^*$ are arbitrary subject to (i)--(iv) above.

\end{note}

\begin{definition}
\label{def:dualtd}
Let $(A,A^*)$ denote a TD pair  on $V$.
Then 
 $(A^*,A)$ is a TD pair
 on $V$, which we refer to
as the  dual of $(A,A^*)$.
\end{definition}
\noindent 
We now consider some examples of TD pairs. 
Our first example is from the theory of association 
schemes.

\begin{example}
\label{ex:pandq} Referring to  
\cite{TersubI},
\cite{TersubII},
\cite{TersubIII},
let $Y$ denote a symmetric 
association scheme, with vertex set $X$.
Assume $Y$ is $P$-polynomial
with respect to an associate matrix $A$, 
and $Q$-polynomial with respect to a primitive
idempotent $E$.
Let 
$\hbox{Mat}_X(\C)$ denote the $\C$-algebra
consisting of all matrices with rows and columns
indexed by $X$, and entries in $\C$.
Fix a vertex $x \in X$, and
let $A^*= A^*(x)$ denote the 
diagonal matrix in 
$\hbox{Mat}_X(\C)$ with  
 diagonal entries
\begin{equation}
A^*_{yy} = 
 \vert X \vert
E_{xy} \qquad \qquad (\forall y \in X).
\end{equation}
Let $T$ denote the subalgebra of 
$\hbox{Mat}_X(\C)$ generated by $A$ and $A^*$, and
let $V$ denote an irreducible $T$-module. Then
 the restrictions  
$A|_V$, $A^*|_V$ form TD pair on $V$. 

\end{example}

\begin{proof} This is immediate from
\cite[ lines (73), (80)]{TersubI}. 

\end{proof}

\noindent  The above example inspired our definition  of
a TD pair, and motivates all that we do in the present paper. 

\medskip
\noindent
Our next example of a TD pair comes from Lie algebra.  

\begin{example}
\label{ex:sl2easy} 
Let $\fld$ denote an algebraically closed field with
characteristic 0, and let $L$ denote the Lie algebra
$sl_2(\fld)$. Let $A$ and $A^*$ denote semi-simple 
elements in $L$, and assume $L$ is generated by
these elements.
Let $V$ denote a finite dimensional
irreducible $L$-module. Then
$A$ and $A^*$ act on $V$ as a TD pair.

\end{example}

\begin{proof} We claim 
 $L$ has a basis $e,f,h$ such that
 $A$ is a nonzero scalar multiple
 of $h$, and such that 
\begin{equation}
\lbrack h,e\rbrack = 2e, \qquad  
\lbrack h,f\rbrack = -2f, \qquad  
\lbrack e,f\rbrack = h,
\label{eq:sl2rel}
\end{equation}
where $\lbrack\,,\,\rbrack $ denotes the Lie bracket.
To see this, we identify $L$ with the Lie algebra
$sl(W)$, where 
$W$ is a two dimensional vector space over $\fld$. Recall
$sl(W)$ is the vector space over $\fld$
consisting
of all linear  transformations
on $W$ that have trace 0, together with  Lie bracket  
$\lbrack x,y\rbrack = xy-yx$. 
Since $A$ is semisimple, 
$W$ has a basis $u,v$ consisting of 
eigenvectors for $A$. Let $r$ denote the eigenvalue
of $A$ associated with $u$.
Since $A$ has trace 0 on $W$, we find 
the eigenvalue of 
$A$ associated with $v$ is $-r$.
 Observe $r\not=0$;
otherwise $A=0$,  contradicting our assumption that $A$ and
$A^*$ generate $L$.
Set $h=r^{-1}A$.
Let $e$ (resp. $f$) denote the linear  transformation
on $W$ satisfying
$ev=u$ and $eu=0$    
(resp. $fu=v$ and $fv=0$).
Observe $e$ and $f$  have trace 0, so they 
are contained in $L$.
 One readily finds $e,f,h$ are linearly independent, and that 
$L$ has dimension 3, so $e,f,h$ form a basis for $L$.
From the construction, we readily obtain 
(\ref{eq:sl2rel}), and 
 our claim is proved.
Using 
(\ref{eq:sl2rel}), we routinely find
the $L$-module $V$ has a basis $v_0, v_1, \ldots, v_d$  
satisfying
\begin{eqnarray}
hv_i &=& (d-2i)v_i \qquad \qquad (0 \leq i \leq d),
\label{eq:honV}
\\
fv_i  &=&(i+1)v_{i+1} \qquad \qquad (0 \leq i < d),\quad  fv_d=0,
\label{eq:fonV}
\\
ev_i  &=& (d-i+1)v_{i-1} \qquad \qquad (0 < i \leq d), \quad  ev_0=0.
\label{eq:eonV}
\end{eqnarray}
From 
(\ref{eq:honV}) we find
$h$ (and $A$) is diagonalizable on $V$.
Let $V_i$ denote the span of $v_i$, for $0\leq i \leq d$.
Then
$V_0, V_1, \ldots, V_d$ is an ordering of the 
maximal eigenspaces of $A$. 
By (\ref{eq:honV})--(\ref{eq:eonV}),
and since 
$A^*$ is a linear combination of $e,f,h$, we
obtain
(\ref{eq:astaractionthreeint}).
We now have Definition 
\ref{def:nonthinlp}(ii).
Interchanging the roles of $A$ and $A^*$ in the above argument,
we find $A^*$ is diagonalizable on $V$, and that 
 Definition 
\ref{def:nonthinlp}(iii) holds.
Since $A, A^*$ generate $L$, and since $V$  is irreducible as
an $L$ module, we obtain Definition 
\ref{def:nonthinlp}(iv).
We have now shown $(A,A^*)$ satisfies the conditions of
Definition \ref{def:nonthinlp}, so
 $(A,A^*)$ is a TD pair on
$V$.

\end{proof}

\noindent  Inspired by the previous example, we 
make the following definition.
Let $(A,A^*)$ denote a 
TD pair.
We say 
$(A,A^*)$ is {\it thin} 
whenever the 
 eigenspaces of both $A$ and $A^*$
all have dimension 1.
In \cite{LS99}, the 
 thin TD pairs are called 
Leonard pairs, and a 
 complete classification of these objects is given.

\medskip
\noindent Here is another example of a TD pair.

\begin{example}
\label{def:onsager}
Let $\fld$ denote an algebraically closed 
field of characteristic 0,
and let $V$ denote a vector space over  $\fld$
with finite positive dimension.
Let $A$ and $A^*$ denote linear
transformations on $V$ satisfying (i)--(ii) below.
\begin{enumerate}
\item There exists  nonzero scalars $b,b^*\in \fld$ such that
\begin{eqnarray}
\lbrack A, \lbrack A, \lbrack A, A^* \rbrack \rbrack \rbrack &=&  b^2\lbrack
A, A^* \rbrack,
\label{eq:DG1}
\\
\lbrack A^*, \lbrack A^*, \lbrack A^*, A \rbrack \rbrack \rbrack &=& 
b^{*2}\lbrack
A^*, A \rbrack,
\label{eq:DG2}
\end{eqnarray}
where 
$\lbrack r,s \rbrack$ means $rs-sr$.
\item There is no subspace $W$ of $V$ such  that  both $AW\subseteq W$,
$A^*W\subseteq W$, other than $W=0$ and $W=V$.
\end{enumerate}
Then $(A,A^*)$ is a TD pair on $V$.
\end{example}

\begin{proof}
For all $\theta\in \fld$,  we define 
\beast
V(\theta) = \lbrace v \in V\;|\;Av = \theta v \rbrace.
\eeast
Observe $V(\theta)$ is nonzero  if and only if $\theta$ is an eigenvalue of
$A$, and in this case 
$V(\theta)$ is the corresponding maximal eigenspace. 
We claim  
\begin{equation}
A^* V(\theta) \subseteq V(\theta-b) + V(\theta) + V(\theta +b)
\qquad \qquad  (\forall \theta \in \fld).
\label{eq:deltaexchange}
\end{equation}
To see 
(\ref{eq:deltaexchange}), let $\theta$ be given, and pick
any $v \in V(\theta)$. 
Expanding 
(\ref{eq:DG1}), we obtain
\begin{equation}
0= A^3A^* - 3A^2A^*A + 3AA^*A^2 - A^*A^3 - b^2 (AA^*-A^*A).
\label{eq:deltadefexpand}
\end{equation}
Applying 
(\ref{eq:deltadefexpand}) to $v$, and using $Av=\theta v$,
we obtain 
\begin{eqnarray}
0 &=& 
 A^3A^*v - 3\theta A^2A^*v + 3\theta^2AA^*v - \theta^3A^*v - 
 b^2 AA^*v+b^2\theta A^*v
\\
&=&
(A-\theta+b )(A-\theta)(A -\theta -b)A^*v.
\label{eq:deltafactor}
\end{eqnarray}
The scalars 
$\theta-b, \theta, \theta+ b$ are distinct, so
apparently
\beast
A^* v \subseteq V(\theta-b) + V(\theta) + V(\theta +b),
\eeast
and (\ref{eq:deltaexchange}) follows.
Since $\fld$ is algebraically closed, all the eigenvalues of
$A$ are contained in $\fld$.
Since $V$ has finite positive
dimension, $A$ has at least one eigenvalue   $\theta$.
Since $\fld$ has characteristic zero,
the scalars  $\theta, \theta-b, \theta -2b, \ldots $
are mutually distinct, so they cannot all be eigenvalues
of $A$; consequently there exists an 
 eigenvalue
$a$ for $A$ such that $a-b$ is not an eigenvalue of $A$. 
Similarily the scalars
$a, a+b, a+2b, \ldots $ are mutually distinct, so they cannot
all be eigenvalues of $A$; consequently 
 there exists a nonnegative integer $d$ such
that $a+bi$ is an eigenvalue of $A$ for $0 \leq i \leq d$, but  
not  for $i=d+1$. Let us abbreviate $V_i = V(a+bi)$
for $0\leq i \leq d$. Then
\begin{equation}
V_0+V_1+\cdots +V_d
\label{eq:sumwillbeV}
\end{equation}
is nonzero and $A$-invariant by construction. 
From
(\ref{eq:deltaexchange}) 
and the construction, we find
\beast
A^* V_i \subseteq V_{i-1} + V_i+ V_{i+1} \qquad \qquad (0 \leq i \leq d),
\eeast
where $V_{-1}=0$, $V_{d+1}=0$.
In particular, 
(\ref{eq:sumwillbeV}) is $A^*$-invariant.
Applying assumption (ii),  we find
the sum
(\ref{eq:sumwillbeV}) equals $V$.
Now apparently 
$A$ is diagonalizable on $V$, and 
Definition
\ref{def:nonthinlp}(ii) holds.
Interchanging the roles of $A$ and $A^*$ in the above argument,
we find $A^*$ is diagonalizable on $V$, 
and 
Definition
\ref{def:nonthinlp}(iii) holds.
Observe 
Definition
\ref{def:nonthinlp}(iv) is just 
assumption (ii) above. 
We have now shown $(A,A^*)$ satisfy the conditions 
of 
Definition
\ref{def:nonthinlp}, so it is a TD pair on $V$.

\end{proof}
\medskip
\noindent  
The equations  
(\ref{eq:DG1}),
(\ref{eq:DG2}) 
are known as the {\it Dolan-Grady}  relations \cite{DateRoan},
\cite{Dav},
\cite{Dolgra}, 
\cite{Roanmpi},
\cite{Ugl} and are important in Statistical Mechanics.
The Lie algebra over $\C$ generated by two symbols $A, A^*$
subject to
(\ref{eq:DG1}),
(\ref{eq:DG2}) 
 (where we interpret 
 $\lbrack \,,\,\rbrack $ as the Lie bracket)
 is infinite dimensional and is known as the 
Onsager algebra 
\cite{CKOns},
\cite{Dav},
\cite{Per},
\cite{Roanmpi},
\cite{Ugl}.

\medskip
\noindent  The next example is a kind of $q$-analog of
the previous one.

\begin{example}
\label{def:serre}
Let $\fld$ denote an algebraically closed 
field, and let $q$
denote a nonzero element in $\fld$ that is not
a root of unity.
Let $V$ denote a vector space over  $\fld$
with finite positive dimension,
and let  $A$ and $A^*$ denote  linear
transformations on $V$ satisfying (i), (ii) below.
\begin{enumerate}
\item $A$ and $A^*$ are not nilpotent, and    
\begin{eqnarray}
0 &=&\lbrack A, A^2A^*-(q+q^{-1})AA^*A+A^*A^2 \rbrack,
\label{eq:Serre1}
\\
0&=&\lbrack A^*, A^{*2}A-(q+q^{-1})A^*AA^*+AA^{*2} \rbrack,
\label{eq:Serre2}
\end{eqnarray}
where $\lbrack r,s\rbrack $ means $rs-sr$. 
\item There is no subspace $W$ of $V$ such  that  both $AW\subseteq W$,
$A^*W\subseteq W$, other than $W=0$ and $W=V$.
\end{enumerate}
Then $(A,A^*)$ is a TD pair on $V$.
%
\end{example}

\begin{proof} The proof is similar to that of 
Example
\ref{def:onsager}.
 For all $\theta\in \fld$,  we define 
\beast
V(\theta) = \lbrace v \in V\;|\;Av = \theta v \rbrace.
\eeast
We claim that for all nonzero $\theta \in \fld$,
\begin{equation}
A^* V(\theta) \subseteq V(q^{-1}\theta) + V(\theta) + V(q\theta).
\label{eq:deltaexchangeserre}
\end{equation}
To see 
(\ref{eq:deltaexchangeserre}), let $\theta$ be given,
and pick any $v \in V(\theta)$.
Expanding the relation (\ref{eq:Serre1}),  
we find
\begin{equation}
0= 
A^3A^*-(q+1+q^{-1})A^2A^*A+(q+1+q^{-1})AA^*A^2-A^*A^3.
\label{eq:deltaserre2}
\end{equation}
Applying  each term in 
(\ref{eq:deltaserre2}) to $v$, 
and using
 $Av=\theta v$,
 we routinely find
\begin{eqnarray}
0 &=& 
 A^3A^*v - (q+1+q^{-1})\theta A^2A^*v + (q+1+q^{-1})\theta^2AA^*v 
 - \theta^3A^*v  
\\
&=&
(A-q^{-1}\theta )(A-\theta )(A -q\theta )A^*v.
\label{eq:deltafactorserre}
\end{eqnarray}
We assume $\theta \not=0$ and that $q$ is not a root of unity,
so 
$q^{-1}\theta, \theta, q\theta $ are distinct. Now  
apparently
\beast
A^* v \subseteq V(q^{-1}\theta) + V(\theta) + V(q\theta),
\eeast
and 
(\ref{eq:deltaexchangeserre}) follows.
Since $\fld$ is algebraically closed, all the eigenvalues of
$A$ are contained in $\fld$.
Since $V$ has finite positive
dimension, and since $A$ is not nilpotent,
$A$ has at least one nonzero eigenvalue   $\theta$.
Since $q$ is not a root of unity,
the scalars  $\theta, q^{-1}\theta, q^{-2}\theta, \ldots $
are mutually distinct, so they cannot all be eigenvalues
of $A$; consequently there exists an 
 eigenvalue
$a$ for $A$ such that $q^{-1}a$ is not an eigenvalue of $A$. 
Similarily the scalars
$a, aq, aq^2, \ldots $ are mutually distinct, so they cannot
all be eigenvalues of $A$; consequently 
 there exists a nonnegative integer $d$ such
that $aq^i$ is an eigenvalue of $A$ for $0 \leq i \leq d$, but  
not  for $i=d+1$. Let us abbreviate $V_i=V(aq^i)$ 
for $0\leq i \leq d$. Then
\begin{equation}
V_0+V_1+\cdots +V_d
\label{eq:sumwillbeVserre}
\end{equation}
is nonzero and $A$-invariant by construction. 
By
(\ref{eq:deltaexchangeserre}) and the construction,
 we find
\beast
A^* V_i \subseteq V_{i-1} + V_i+ V_{i+1} \qquad \qquad (0 \leq i \leq d),
\eeast
where $V_{-1}=0$, $V_{d+1}=0$.
In particular, 
(\ref{eq:sumwillbeVserre})
is $A^*$-invariant.
Applying assumption (ii) above,  we find
the sum
(\ref{eq:sumwillbeVserre})
equals $V$.
Now apparently 
$A$ is diagonalizable on $V$, and 
Definition
\ref{def:nonthinlp}(ii) holds.
Interchanging the roles of $A$ and $A^*$ in the above argument,
we find $A^*$ is diagonalizable on $V$, 
and 
Definition
\ref{def:nonthinlp}(iii) holds.
Observe 
Definition
\ref{def:nonthinlp}(iv) is just 
assumption (ii) above. 
We have now shown $(A,A^*)$ satisfy the conditions 
of 
Definition
\ref{def:nonthinlp}, so it is a TD pair on $V$.

\end{proof}

\medskip
\noindent
The relations 
(\ref{eq:Serre1}),
(\ref{eq:Serre2})
are known as {\it quantum Serre relations}.
They are among the defining relations for 
the quantum affine algebra $U_r({\widehat {sl}}_2)$, where
$r^2=q$. See 
\cite{CPqaa}, \cite{CPqaar} for more information about  this algebra.


\section{TD systems}
\noindent 
When working with a TD pair, it is often  convenient to consider 
a closely  related but somewhat more abstract concept, called a TD system.
To define it, we recall a few concepts from linear algebra.
Let $V$ denote a  vector space over
$\fld$ with finite positive dimension, and
let $\hbox{End}(V)$ denote the $\fld$-algebra consisting
of all linear transformations on $V$.
let $A$ denote a  diagonalizable element in 
$\hbox{End}(V)$.
Let $\theta_0, \theta_1, \ldots, \theta_d$ denote an ordering of 
the eigenvalues
of $A$, and put 
\beast
E_i = \prod_{{0 \leq  j \leq d}\atop
{j\not=i}} {{A-\theta_j I}\over {\theta_i-\theta_j}}
\eeast
for $0 \leq i \leq d$.
By elementary linear algebra,
\begin{eqnarray}
&&AE_i = E_iA = \theta_iE_i \qquad \qquad  (0 \leq i \leq d),
\label{eq:primid1S99}
\\
&&
\quad E_iE_j = \delta_{ij}E_i \qquad \qquad (0 \leq i,j\leq d),
\label{eq:primid2S99}
\\
&&
\qquad \qquad \sum_{i=0}^d E_i = I.
\label{eq:primid3S99}
\end{eqnarray}
From this, 
one  finds  $E_0, E_1, \ldots, E_d$ is a  basis for the
subalgebra of $\hbox{End}(V)$ generated by $A$.
We refer to $E_i$ as the {\it primitive idempotent} of
$A$ associated with $\theta_i$.
It is helpful to think of these primitive idempotents as follows. 
From 
(\ref{eq:primid2S99}), 
(\ref{eq:primid3S99}) one readily finds 
\begin{eqnarray}
V = E_0V + E_1V + \cdots + E_dV \qquad \qquad (\hbox{direct sum}).
\label{eq:VdecompS99}
\end{eqnarray}
For $0\leq i \leq d$, $E_iV$ is the  eigenspace of
$A$ associated with the 
eigenvalue $\theta_i$, 
and $E_i$ is the projection of $V$  onto this eigenspace.

\begin{definition}
\label{def:deflstalkS99}
Let $\fld $ denote a field, and let $V$ denote a vector space
over $\fld$ with finite positive dimension. 
By a {\it Tridiagonal system} (or TD system) on $V$,  we mean a 
sequence 
\begin{equation}
\; \Phi = (A;\,E_0,\,E_1,\,\ldots,
\,E_d;\,A^*;\,E^*_0,\,E^*_1,\,\ldots,\,E^*_\delta)
\label{eq:ourstartingpt}
\end{equation}
 that satisfies  (i)--(vi) below.
\begin{enumerate}
\item $A$,  $\;A^*\;$ are both diagonalizable  linear transformations on $V$. 
\item $E_0,\,E_1,\,\ldots,\,E_d\;$ is an ordering of the primitive 
idempotents of $\;A$.
\item $E^*_0,\,E^*_1,\,\ldots,\,E^*_{\delta}\;$ is an ordering of the primitive 
idempotents of $\;A^*$.
\item 
$E_iA^*E_j = 0 \quad \hbox{if}\quad |i-j|>1, \qquad  \qquad    
(0 \leq i,j\leq d)$.
\item 
 $E^*_iAE^*_j = 0 \quad \hbox{if}\quad \vert i-j\vert > 1,
\qquad \qquad 
(0 \leq i,j\leq \delta).$
\item There is no subspace $W$ of $V$ such  that  both $AW\subseteq W$,
$A^*W\subseteq W$, other than $W=0$ and $W=V$.
\end{enumerate}
We refer to $d$ as the {\it diameter} of $\Phi$,  
and say $\Phi$ is {\it over } $\fld$. 
 For notational convenience, we set $E_{-1}=0$, $E_{d+1}=0$, $
E^*_{-1}=0$, $E^*_{\delta+1}=0$.

\end{definition}

\noindent Referring to Definition 
\ref{def:deflstalkS99}, we 
do not assume the primitive idempotents of $A$ and $A^*$ all 
have rank 1.  
A TD system for which   these primitive idempotents all have rank 
1 is called a Leonard system \cite{LS99}.   The Leonard systems  are completely classified by
Terwilliger \cite{LS99}.

\medskip
\noindent In the next two lemmas, we give the connection between
TD pairs and TD systems.

\begin{lemma}
\label{lem:tdpairsys}
Let $\fld $ denote a field, and let $V$ denote a vector space
over $\fld$ with finite positive dimension. 
Let $(A,A^*)$ denote a TD pair on $V$.  Let $V_0, V_1,\ldots ,V_d$
denote an ordering of the eigenspaces of $A$ satisfying
(\ref{eq:astaractionthreeint}), and for $0 \leq i \leq d$ let
$E_i$ denote the projection of $V$ onto $V_i$. 
 Let $V^*_0, V^*_1,\ldots ,V^*_\delta $
denote an ordering of the eigenspaces of $A^*$ satisfying
(\ref{eq:aactionthreeint}),
and for $0 \leq i \leq \delta $ let
$E^*_i$ denote the projection of $V$ onto $V^*_i$. 
Then
\beast
(A;E_0,E_1,\ldots,E_d;A^*;E^*_0,E^*_1,\ldots,E^*_{\delta})
\eeast
 is a  TD system on $V$.

\end{lemma}

\begin{proof} We verify the conditions (i)--(vi) of Definition 
\ref{def:deflstalkS99}. Condition (i) is just
 Definition 
\ref{def:nonthinlp}(i).
Conditions (ii) and (iii)  follow from the comment after 
(\ref{eq:VdecompS99}).
To obtain condition (iv), pick any integers $i,j$
$(0 \leq i,j\leq d)$, and assume  $\vert i-j \vert >1$.  For all
$v \in V$, 
\beast
E_iA^*E_jv &\in &  E_iA^*V_j
\\
&\subseteq &  E_i(V_{j-1} + V_j + V_{j+1})
\\
&=& 0,
\eeast
so $E_iA^*E_j=0$.  We have now verified condition (iv), and
condition (v) is similarily obtained.
Condition (vi) is just Definition 
\ref{def:nonthinlp}(iv).

\end{proof}

\begin{lemma}
\label{lem:tdsyspair}
Let $\fld $ denote a field, let $V$ denote a vector space
over $\fld$ with finite positive dimension, and
let 
\beast
\Phi = (A;E_0,E_1,\ldots,E_d;A^*;E^*_0,E^*_1,\ldots,E^*_{\delta})
\eeast
denote a TD sytem on $V$. Then
\begin{eqnarray}
&&A^*E_iV \subseteq E_{i-1}V + E_iV+E_{i+1}V \qquad \qquad (0 \leq i \leq d),
\label{eq:asspread}
\\
&&AE^*_iV \subseteq E^*_{i-1}V + E^*_iV+E^*_{i+1}V 
\qquad \qquad (0 \leq i \leq \delta).
\label{eq:aspread}
\end{eqnarray}
Moreover, $(A,A^*)$ is a TD pair on $V$. We refer to $(A,A^*)$
as the TD pair associated with $\Phi$. 
\end{lemma}

\begin{proof} To obtain
(\ref{eq:asspread}), observe by 
(\ref{eq:primid3S99}) and 
Definition \ref{def:deflstalkS99}(iv) that 
\beast
A^*E_iV &=& (E_0+E_1+\cdots + E_d)A^*E_iV
\\
&=& 
(E_{i-1}+E_i+E_{i+1})A^*E_iV
\\
&\subseteq & E_{i-1}V + E_iV+E_{i+1}V.
\eeast
Line (\ref{eq:aspread}) is similarily obtained, and the last
assertion follows. 

\end{proof}

\noindent We finish this section with an observation.

\begin{lemma}
\label{lem:eaenotzero}
Referring to the TD system 
(\ref{eq:ourstartingpt}), we have
\begin{enumerate}
\item 
$E_iA^*E_j \not= 0 \quad \hbox{if}\quad |i-j|=1, \qquad  \qquad    
(0 \leq i,j\leq d)$,
\item 
 $E^*_iAE^*_j \not= 0 \quad \hbox{if}\quad \vert i-j\vert = 1,
\qquad \qquad 
(0 \leq i,j\leq \delta).$
\end{enumerate}
\end{lemma}

\begin{proof} (i) Pick any integer $i$ $(1 \leq i \leq d)$, and first assume
$E_{i}A^*E_{i-1}=0$. Put $W:=E_0V+E_{1}V + \cdots + E_{i-1}V$.
By 
(\ref{eq:VdecompS99}), and since each of
$E_0V, E_1V, \ldots, E_dV$
is nonzero, we find 
$W\not=0$, $W\not=V$.   
Observe 
$AW\subseteq W$ by the construction.
Combining 
(\ref{eq:asspread})
and our assumption, we find
$A^*W \subseteq W$,
and now condition (vi) of Definition
\ref{def:deflstalkS99} 
is contradicted. Next suppose  
$E_{i-1}A^*E_{i}=0$. In this case we routinely find
$W:=E_iV+E_{i+1}V + \cdots + E_dV$ provides
a contradiction to 
condition (vi) of Definition
\ref{def:deflstalkS99}.

\noindent (ii) Similar to the proof of (i) above.

\end{proof}

\section{The relatives of a TD system} 

\medskip
\noindent A given TD system  can be modified in  several
ways to get a new TD system. For instance, 
let $\Phi$ 
 denote the TD system in
(\ref{eq:ourstartingpt}), 
 and let $\alpha$, $\alpha^*$, $\beta$,
$\beta^*$ denote scalars in $\fld$ such that $\alpha \not=0$, $\alpha^*\not=0$.
Then
\beast
(\alpha A+\beta I; E_0, E_1, \ldots, E_d; 
\alpha^* A^*+\beta^* I; E^*_0, E^*_1, \ldots, E^*_{\delta})
\eeast
is a TD system on $V$.
Also, 
\begin{eqnarray}
 \;\Phi^*&:=& (A^*; E^*_0,E^*_1,\ldots,E^*_{\delta};A;E_0,E_1, \ldots,E_d),
\label{eq:lsdualS99}
\\
\Phi^{\downarrow}&:=& (A; E_0,E_1,\ldots,E_d;A^*;E^*_{\delta},E^*_{\delta-1}, \ldots,E^*_0),
\label{eq:lsinvertS99}
\\
\Phi^{\Downarrow} 
&:=& (A; E_d,E_{d-1},\ldots,E_0;A^*;E^*_0,E^*_1, \ldots,E^*_{\delta})
\label{eq:lsdualinvertS99}
\end{eqnarray}
are TD  systems on $V$.
 We refer to $\Phi^*$
(resp.  
 $\Phi^\downarrow$)   
(resp.  
 $\Phi^\Downarrow$) 
 as the  
{\it dual} 
(resp. {\it first inversion})
(resp.  {\it second inversion}) of  $\Phi$.
Viewing $*, \downarrow, \Downarrow$
as permutations on the set of all TD systems,
\begin{eqnarray}
&&\qquad \qquad \qquad  *^2 \;=\;  
\downarrow^2\;= \;
\Downarrow^2 \;=\;1,
\qquad \quad 
\label{eq:deightrelationsAS99}
\\
&&\Downarrow *\; 
=\;
* \downarrow,\qquad \qquad   
\downarrow *\; 
=\;
* \Downarrow,\qquad \qquad   
\downarrow \Downarrow \; = \;
\Downarrow \downarrow.
\qquad \quad 
\label{eq:deightrelationsBS99}
\end{eqnarray}
The group generated by symbols 
$*, \downarrow, \Downarrow $ subject to the relations
(\ref{eq:deightrelationsAS99}),
(\ref{eq:deightrelationsBS99})
is the dihedral group $D_4$.  
We recall $D_4$ is the group of symmetries of a square,
and has 8 elements.
Apparently $*, \downarrow, \Downarrow $ induce an action of 
 $D_4$ on the set of all TD systems.
Two TD systems will be called {\it relatives} whenever they
are in the same orbit of this $D_4$ action.
The relatives of $\Phi$ are as follows:
\medskip

\centerline{
\begin{tabular}[t]{c|c}
        name &relative \\ \hline 
        $\Phi$ & $(A;E_0,E_1,\ldots,E_d;A^*;E^*_0,E^*_1,\ldots,E^*_\delta)$   \\ 
        $\Phi^\downarrow$ &
         $(A;E_0,E_1,\ldots,E_d;A^*;E^*_\delta,E^*_{\delta-1},\ldots,E^*_0)$   \\ 
        $\Phi^\Downarrow$ &
         $(A;E_d,E_{d-1},\ldots,E_0;A^*;E^*_0,E^*_1,\ldots,E^*_\delta)$   \\ 
        $\Phi^{\downarrow \Downarrow}$ &
         $(A;E_d,E_{d-1},\ldots,E_0;A^*;E^*_\delta,E^*_{\delta-1},\ldots,E^*_0)$   \\ 
	$\Phi^*$ &
        $(A^*;E^*_0,E^*_1,\ldots,E^*_\delta;A;E_0,E_1,\ldots,E_d)$   \\ 
        $\Phi^{\downarrow *}$ &
	 $(A^*;E^*_\delta,E^*_{\delta-1},\ldots,E^*_0;  
         A;E_0,E_1,\ldots,E_d)$ \\
        $\Phi^{\Downarrow *}$ &
	 $(A^*;E^*_0,E^*_1,\ldots,E^*_\delta;    
         A;E_d,E_{d-1},\ldots,E_0)$ \\ 
	$\Phi^{\downarrow \Downarrow *}$ &
	 $(A^*;E^*_\delta,E^*_{\delta-1},\ldots,E^*_0;    
         A;E_d,E_{d-1},\ldots,E_0)$
	\end{tabular}}
\medskip
\noindent 

%
%
%

\medskip
\noindent   With reference to 
Definition \ref{def:dualtd} and 
Lemma \ref{lem:tdsyspair},
it is routine to show two TD systems are relatives
if and only if their associated TD pairs are equal
or dual.

\noindent

\medskip 
\noindent We now introduce two sequences of parameters
that we will use to describe a given TD system.

\begin{definition}
\label{def:lsdefcommentsS99}
Let $\Phi$ denote the TD system in 
(\ref{eq:ourstartingpt}).
For $0 \leq i \leq d$, 
we let $\theta_i $ 
denote the eigenvalue
of $A$
associated with $E_i$.
We refer to  $\theta_0, \theta_1, \ldots, \theta_d$ as the 
eigenvalue sequence of $\Phi$.
For $0 \leq i \leq \delta$, 
we let $\theta^*_i $ 
denote the eigenvalue
of $A^*$
associated with $E^*_i$.
We refer to  $\theta^*_0, \theta^*_1, \ldots, \theta^*_\delta$ as the 
dual eigenvalue sequence of $\Phi$.
We remark 
 $\theta_0, \theta_1, \ldots, \theta_d$ are mutually distinct, and 
$\theta^*_0, \theta^*_1, \ldots, \theta^*_\delta$ are mutually distinct. 

\end{definition}
\medskip

\section{\bf The split decomposition} 


\begin{definition}
\label{def:splitconsetupS99}
In this section, $\fld $ will 
denote a field, and $V$ will denote a  vector space over $\fld$ with finite
positive dimension.
We let 
\beast
\Phi =(A;E_0,E_1,\ldots, E_d;A^*;E^*_0,E^*_1,\ldots,E^*_{\delta})
\eeast
denote a TD system on $V$, with eigenvalue  sequence
$\theta_0, \theta_1, \ldots, \theta_d$ and dual eigenvalue 
sequence
 $\theta^*_0, \theta^*_1, \ldots, \theta^*_{\delta}$.
\end{definition}

\noindent Referring to the above definition, we will show $\delta = d$. We will then
show 
there exists 
a unique sequence $U_0, U_1, \ldots, U_d $ of  subspaces of $V$ such that 
\begin{eqnarray}
&&\qquad V = U_0+ U_1 + \cdots + U_d \qquad \qquad (\hbox{direct sum}),
\label{eq:vdecomp}
\\
&&(A-\theta_iI)U_i\subseteq U_{i+1} \qquad \quad (0 \leq i <d),
\qquad (A-\theta_dI)U_d=0,   
\label{eq:Aaction}
\\
&& (A^*-\theta^*_iI)U_i\subseteq U_{i-1} \qquad \quad (0 < i \leq d), 
\qquad (A^*-\theta^*_0I)U_0=0.    
\label{eq:Asaction}
\end{eqnarray}

\noindent The following notation will be useful.
%
%

\begin{definition}
\label{def:VijS99}
With reference to 
Definition
\ref{def:splitconsetupS99}, 
we set 
\begin{equation}
V_{ij} =  \Biggl(\sum_{h=0}^i E^*_hV\Biggr)\cap
 \Biggl(\sum_{k=j}^d E_kV\Biggr)
\label{eq:defofvijS99}
\end{equation}
for all integers $i,j$. We interpret the sum on the left in
(\ref{eq:defofvijS99}) to be 0 (resp. $V$) if $i<0$  (resp. $i>\delta $).
Similarily, we interpret the sum on the right in
(\ref{eq:defofvijS99}) to be V (resp. $0$) if $j<0$  (resp. $j>d$).

\end{definition}

\begin{lemma}
\label{lem:thevijbasicfactsS99}
With reference to 
Definition
\ref{def:splitconsetupS99} and 
Definition \ref{def:VijS99}, we have
\begin{enumerate}
\item
$V_{i0} = 
 E^*_0V+E^*_1V+\cdots+E^*_iV \qquad \qquad  (0 \leq i \leq \delta),$
\item
$V_{\delta j} = 
 E_jV+E_{j+1}V+\cdots+E_dV \qquad \qquad (0 \leq j\leq d).$
\end{enumerate}

\end{lemma}

\begin{proof} 
To get (i), set $j=0$ in 
(\ref{eq:defofvijS99}), and apply 
(\ref{eq:VdecompS99}). Line  (ii) is similarily obtained.
\end{proof}

\begin{lemma}
\label{lem:howaactsonvijS99}
With reference to 
Definition
\ref{def:splitconsetupS99} 
and Definition \ref{def:VijS99}, the following (i)--(iv) hold
for $0 \leq i\leq \delta$ and $0 \leq j\leq d$.
\begin{enumerate}
\item
$(A-\theta_jI)V_{ij} \subseteq V_{i+1,j+1}$,
\item
$AV_{ij} \subseteq V_{ij} + V_{i+1,j+1}$, 
\item
$(A^*-\theta^*_iI)V_{ij} \subseteq V_{i-1,j-1}$,
\item
$A^*V_{ij} \subseteq V_{ij} + V_{i-1,j-1}$.
\end{enumerate}
\end{lemma}

\begin{proof} (i) Using 
(\ref{eq:aspread}),
\begin{equation}
(A-\theta_jI)\sum_{h=0}^i E^*_hV \subseteq 
\sum_{h=0}^{i+1} E^*_hV,
\label{eq:aminusthetajAS99}
\end{equation}
and using
(\ref{eq:primid1S99}), 
\begin{equation}
(A-\theta_jI)\sum_{k=j}^d E_kV = 
\sum_{k=j+1}^d E_kV.
\label{eq:aminusthetajBS99}
\end{equation}
Evaluating $(A-\theta_jI)V_{ij}$ using 
(\ref{eq:defofvijS99}),
(\ref{eq:aminusthetajAS99}), 
(\ref{eq:aminusthetajBS99}), 
 we
routinely find it is contained in $V_{i+1,j+1}$.

\noindent (ii) Immediate from  (i) above.

\noindent (iii)
Using (\ref{eq:primid1S99}),
\begin{equation}
(A^*-\theta^*_iI)\sum_{h=0}^i E^*_hV = 
\sum_{h=0}^{i-1} E^*_hV,
\label{eq:aminusthetajAS99ag}
\end{equation}
and using
(\ref{eq:asspread}),
\begin{equation}
(A^*-\theta^*_iI)\sum_{k=j}^d E_kV \subseteq 
\sum_{k=j-1}^d E_kV.
\label{eq:aminusthetajBS99ag}
\end{equation}
Evaluating $(A^*-\theta^*_iI)V_{ij}$ using 
(\ref{eq:defofvijS99}),
(\ref{eq:aminusthetajAS99ag}), 
(\ref{eq:aminusthetajBS99ag}), 
 we
routinely find it is contained in $V_{i-1,j-1}$.

\noindent (iv) Immediate from (iii) above.

\end{proof}

\begin{lemma}
\label{lem:manyvijzeroS99}
The scalars 
 $d$ and $\delta $ 
from Definition
\ref{def:splitconsetupS99}
 are equal. Moreover,
with reference to 
Definition \ref{def:VijS99}, 
\begin{equation}
V_{ij}= 0 \quad \hbox{if}\quad i<j,\qquad  \qquad (0 \leq i,j\leq d).
\label{eq:manyvijzeroS99}
\end{equation}
\end{lemma}

\begin{proof} 
Switching $A$ and $A^*$ if necessary, we may assume
$\delta \leq d$. We first show 
(\ref{eq:manyvijzeroS99}). To do this, 
we show the sum 
\begin{equation}
V_{0r}+V_{1,r+1}+\cdots +V_{d-r,d}
\label{eq:manyvijzeroAS99}
\end{equation}
is zero  for $0 <r \leq d$. Let $r$ be given, and let $W$ denote the sum
in 
(\ref{eq:manyvijzeroAS99}). Applying Lemma 
\ref{lem:howaactsonvijS99}(ii),(iv), we find
$AW\subseteq W$ and $A^*W\subseteq W$.
By the definition of a TD system, we find
$W=0$ or $W=V$.
 By Definition
\ref{def:VijS99}, each term in 
(\ref{eq:manyvijzeroAS99}) is contained in 
\begin{equation}
E_rV+E_{r+1}V+\cdots +E_dV,
\label{eq:manyvijzeroBS99}
\end{equation}
so $W$ is contained in 
(\ref{eq:manyvijzeroBS99}). The sum 
(\ref{eq:manyvijzeroBS99}) is properly contained in $V$ by 
(\ref{eq:VdecompS99}), and since $r>0$.
Apparently $W\not=V$, so
$W=0$. We have now shown  
(\ref{eq:manyvijzeroAS99}) is zero for $0<r\leq d$, and 
(\ref{eq:manyvijzeroS99}) follows.
It remains to show $\delta = d$.
Suppose $\delta \not=d$, so that $\delta < d$ by our initial assumption.
On one hand, 
setting $i=\delta$, $j=d$ in 
(\ref{eq:manyvijzeroS99}), we find $V_{\delta d} =0$. On the other hand 
$V_{\delta d} =E_dV$ by Lemma 
\ref{lem:thevijbasicfactsS99}(ii), so $V_{\delta d} \not=0$.
We now have a contradiction,
so $\delta = d$.

\end{proof}

\begin{theorem}
\label{thm:maincharls}
With reference to 
Definition
\ref{def:splitconsetupS99}, let 
$U_0, U_1,\ldots, U_d$ denote any subspaces of $V$. Then the following are equivalent.
\begin{enumerate}
\item $U_i = (E^*_0V+E^*_1V+\cdots + E^*_iV)\cap (E_iV+E_{i+1}V+\cdots + E_dV) \qquad (0 \leq i \leq d)$.
\item  The sequence $U_0, U_1, \ldots, U_d $ satisfies 
(\ref{eq:vdecomp}),
(\ref{eq:Aaction}), 
(\ref{eq:Asaction}).
\item For $0 \leq i \leq d,$ both 
\begin{eqnarray}
U_i+U_{i+1}+\cdots +U_d &=& E_iV+E_{i+1}V +\cdots + E_dV,
\label{eq:vsumitod}
\\
U_0+U_1+\cdots +U_i &=& E^*_0V+E^*_1V +\cdots + E^*_iV. 
\label{eq:vsumzeroi}
\end{eqnarray}
\end{enumerate}
\end{theorem}

\begin{proof} $(i)\rightarrow (ii)$  
To get 
(\ref{eq:Aaction}) and 
(\ref{eq:Asaction}), set $j=i$ in Lemma 
\ref{lem:howaactsonvijS99}(i),(iii),  and observe
$U_i = V_{ii}$.
To obtain
(\ref{eq:vdecomp}),
let $W$ denote the sum on the right in 
that line.
Then $AW\subseteq W$ by 
(\ref{eq:Aaction}), 
and $A^*W\subseteq W$ by 
(\ref{eq:Asaction}). 
By the definition of a TD system,
we find
$W=0$ or $W=V$.  $W$ contains 
 $U_{0}$, and $U_{0}=E^*_0V$ is nonzero,
 so
$W\not=0$.
It follows $W=V$, and in other words 
\begin{equation}
V= U_{0}+U_{1}+\cdots+U_{d}.
\label{eq:vijgivesdirsumCS99}
\end{equation}
We show the sum 
(\ref{eq:vijgivesdirsumCS99}) is direct.  To do this,
we show
\beast
 (U_{0}+U_{1}+\cdots +U_{i-1})\cap U_{i} =0
\eeast
for $1 \leq i \leq d$.
Let the integer $i$ be given. From the construction 
\beast
U_{j} \subseteq E^*_0V+E^*_1V+\cdots +E^*_{i-1}V
\eeast
for $0 \leq j \leq i-1$, and
\beast
U_{i}\subseteq E_iV+E_{i+1}V+\cdots +E_dV.
\eeast
It follows
\beast
 &&(U_{0}+U_{1}+\cdots +U_{i-1})\cap U_{i}
\\
 && \qquad  \subseteq
 (E^*_0V+E^*_1V+\cdots +E^*_{i-1}V)\cap
(E_iV+E_{i+1}V+\cdots +E_dV) \qquad \qquad 
\\
&& \qquad  = V_{i-1,i}
\\
&& \qquad  = 0
\eeast
in view of Lemma 
\ref{lem:manyvijzeroS99}.
We have now shown
the sum (\ref{eq:vijgivesdirsumCS99}) is direct,
so 
(\ref{eq:vdecomp}) holds.

\noindent 
 $(ii)\rightarrow (iii)$  
First consider 
(\ref{eq:vsumitod}).  
Let $i$ be given, and abbreviate
\beast
Z&=& E_iV+E_{i+1}V+\cdots + E_dV, \qquad \quad 
W=U_i+U_{i+1}+\cdots + U_d.
\eeast
We show $Z=W$. To obtain $Z\subseteq W$, set 
$X=\prod_{h=0}^{i-1} (A-\theta_hI)$, and observe
$Z=XV$ by 
(\ref{eq:primid1S99})--(\ref{eq:primid3S99}). Using 
(\ref{eq:Aaction}), we find $XU_j\subseteq W$ for $0 \leq j \leq d$, so $XV\subseteq W$ in view of
(\ref{eq:vdecomp}).
We now have 
$Z\subseteq W$. To obtain $W\subseteq Z$, set 
$Y=\prod_{h=i}^{d} (A-\theta_hI)$, and observe
\begin{eqnarray}
Z&=&\lbrace v \in V \;|\;Yv = 0\rbrace .
\label{eq:Keraction}
\end{eqnarray}
Using 
(\ref{eq:Aaction}), we find 
$YU_j=0$ for $i \leq j \leq d$, so $YW=0$. Combining 
this with
(\ref{eq:Keraction}), we find $W\subseteq Z$. We now have $Z=W$ and hence
(\ref{eq:vsumitod}).
Line (\ref{eq:vsumzeroi})  is similarily obtained.

\noindent $(iii)\rightarrow (i)$ We first
show $U_0, U_1, \ldots, U_d$ are linearly
independent. To do this, we show
\begin{equation} 
 (U_{0}+U_{1}+\cdots +U_{i-1})\cap U_{i}
\label{eq:linindepag}
\end{equation}
is zero for $1 \leq i \leq d$. Let $i$ be given. From
(\ref{eq:vsumitod}), 
 (\ref{eq:vsumzeroi}),
we find (\ref{eq:linindepag}) is contained in 
\begin{equation}
 (E^*_0V+E^*_1V+\cdots +E^*_{i-1}V)\cap
(E_iV+E_{i+1}V+\cdots +E_dV). \qquad \qquad 
\label{eq:thatintagain}
\end{equation}
The expression (\ref{eq:thatintagain}) equals $V_{i-1,i}$, and is hence zero
by Lemma 
\ref{lem:manyvijzeroS99}. It follows 
(\ref{eq:linindepag}) is zero, and we have now shown
 $U_0, U_1, \ldots, U_d$ are linearly independent.
Combining  this with  
(\ref{eq:vsumitod}), 
 (\ref{eq:vsumzeroi}), we find
\beast
U_i &=& (U_0+U_1+\cdots + U_i)\cap (U_i+U_{i+1}+\cdots + U_d)
\\
&=& 
 (E^*_0V+E^*_1V+\cdots +E^*_iV)\cap
(E_iV+E_{i+1}V+\cdots +E_dV),
\eeast
as desired.

\end{proof}

\section{\bf Some projections}

\begin{definition}
\label{def:meaningofVi}
In this section, $\fld $ will 
denote a field, and $V$ will denote a  vector space over $\fld$ with finite
positive dimension.
We let 
\begin{equation}
\Phi =(A;E_0,E_1,\ldots, E_d;A^*;E^*_0,E^*_1,\ldots,E^*_d)
\label{eq:phidisplay}
\end{equation}
denote a TD system on $V$, with eigenvalue  sequence
$\theta_0, \theta_1, \ldots, \theta_d$ and dual eigenvalue 
sequence
 $\theta^*_0, \theta^*_1, \ldots, \theta^*_d$.
We let $U_0, U_1, \ldots, U_d$ denote the subspaces of $V$ that
satisfy (i)--(iii) of Theorem
\ref{thm:maincharls}.
\end{definition}

\begin{definition}
\label{def:deffi}
With reference to 
Definition
\ref{def:meaningofVi},
for $0 \leq i \leq d$, we let 
$F_i$ denote the linear transformation on $V$
satisfying both
\begin{eqnarray}
&&\qquad (F_i-I)U_i=0,
\label{eq:fiui}
\\
&&F_iU_{j}=0 \quad \hbox{if}\quad j\not=i, \qquad (0 \leq j \leq d).
\label{eq:fiuj}
\end{eqnarray}
In other words, $F_i$  is the projection map from $V$ onto $U_i$.
For notational convenience, we define $F_{-1} = 0 $ and $F_{d+1} = 0$.
\end{definition}

\begin{lemma}
\label{lem:fielem}
With reference to Definition
\ref{def:meaningofVi}
and Definition \ref{def:deffi},
\begin{eqnarray}
&&F_iF_j = \delta_{ij}F_i \qquad \qquad (0 \leq i,j\leq d),
\label{eq:fifiidem}
\\
&&\qquad F_0+F_1+\cdots +F_d = I,
\label{eq:fisumtoident}
\\
&&F_iV = U_{i} \qquad \qquad (0 \leq i \leq d).
\label{eq:fivjustvi}
\end{eqnarray}

\end{lemma}

\begin{proof} Immediate from Definition
\ref{def:deffi} and 
(\ref{eq:vdecomp}).

\end{proof}

\begin{lemma}
\label{def:fielemprop}
With reference to 
Definition
\ref{def:meaningofVi}
and Definition \ref{def:deffi},
the products
$E_iF_j$, $F_iE_j$, $E^*_jF_i$, $F_jE^*_i$ 
are all zero
for $0 \leq i<j\leq d$.

\end{lemma}
\begin{proof} Let the integers $i,j$ be given. Combining  
(\ref{eq:primid2S99}),
(\ref{eq:vsumitod}),
(\ref{eq:fivjustvi}), 
 we find
\beast
E_i F_jV&=&E_iU_j         
\\
&\subseteq& E_i(E_jV+E_{j+1}V + \cdots + E_dV) 
\\
&=& 0, 
\eeast
so $E_iF_j=0$. Similarily, using
(\ref{eq:vsumitod}) and 
(\ref{eq:fiuj}),
\beast
F_iE_jV &\subseteq & F_i(E_jV+E_{j+1}V + \cdots + E_dV)
\\
&= & F_i(U_j+U_{j+1} + \cdots + U_d)
\\
&=& 0,
\eeast
so $F_iE_j=0$. The remaining assertions are similarily proved.

\end{proof}

\begin{lemma}
\label{def:twoisoprel}
With reference to 
Definition
\ref{def:meaningofVi}
and Definition
\ref{def:deffi},
we have
\begin{eqnarray}
F_iE_iF_i = F_i, \qquad \qquad E_iF_iE_i  = E_i,
\label{eq:twoisoprelB}
\\
F_iE^*_iF_i = F_i, \qquad \qquad E^*_iF_iE^*_i  = E^*_i,
\label{eq:twoisoprelA}
\end{eqnarray}
for $0 \leq i \leq d$. 
\end{lemma}

\begin{proof} To obtain the equation on the left in
(\ref{eq:twoisoprelB}), 
we evaluate the expression
\begin{equation}
F_i(E_0+E_1+\cdots + E_d)F_i
\label{eq:fisumfi}
\end{equation}
in two ways.
First, the expression
(\ref{eq:fisumfi}) equals
 $F_i$,
since the sum in the middle is the identity, and since $F_i^2=F_i$.
Pick any integer $j$ $(0 \leq j \leq d)$. 
By Lemma 
\ref{def:fielemprop},
we have $F_iE_j= 0$ if 
$j>i$, and 
$E_jF_i = 0 $ if $j<i$, so 
$F_iE_jF_i= 0 $ if $j\not=i$. 
Apparently 
(\ref{eq:fisumfi}) equals $F_iE_iF_i$, and we now 
have the 
equation on the left in 
(\ref{eq:twoisoprelB}). The  remaining  
equations are proved in a similar manner.

\end{proof}

\noindent We interpret the above lemma as follows.

\begin{lemma}
\label{def:efisofin}
With reference to 
Definition
\ref{def:meaningofVi}
and Definition
\ref{def:deffi}, the following (i), (ii) hold  for $0 \leq i \leq d$. 
\begin{enumerate}
\item The linear transformations
\beast
{{U_i\quad \rightarrow \quad  E_iV}\atop {v \quad \rightarrow \quad E_iv}}
\qquad \qquad \qquad 
{{E_iV\quad \rightarrow \quad U_i}\atop {v \quad \; \rightarrow \quad F_iv}}
\eeast
are bijections, and  moreover, they are inverses.
\item The linear transformations 
\beast
{{U_i\quad \rightarrow \quad  E^*_iV}\atop {v \quad \rightarrow \quad E^*_iv}}
\qquad \qquad \qquad 
{{E^*_iV\quad \rightarrow \quad U_i}\atop {v \quad \; \rightarrow \quad F_iv}}
\eeast
are bijections, and  moreover, they are inverses.
\end{enumerate}
\end{lemma}

\begin{proof} (i) 
The maps are 
inverses in view of 
(\ref{eq:twoisoprelB}). It follows they are bijections. 

\noindent (ii) Similar to the proof of (i).

\end{proof}

\begin{corollary}
\label{thm:symunimodal} 
With reference to Definition
\ref{def:meaningofVi},
for $0 \leq i \leq d$, the dimensions of 
$E_iV$, $U_i$, and $E^*_iV$ are equal. Denoting
this common dimension  by $\rho_i$, we have 
$\rho_i = \rho_{d-i}$.
\end{corollary}

\begin{proof} It is immediate from Lemma
\ref{def:efisofin} that the dimensions of
$E_iV, U_i, E^*_iV$ are equal. 
Denote this common dimension by $\rho_i$.
To  show $\rho_i=\rho_{d-i}$, we show 
 $E_iV$ and $E^*_{d-i}V$ have the same dimension.
We just showed $E_iV$ and $E^*_iV$ have the same dimension.
Applying this result to $\Phi^\downarrow$, we find
$E_iV$ and $E^*_{d-i}V$ have the same dimension. 

\end{proof}

\noindent We finish this section with a remark.

\begin{lemma}
\label{lem:thevifinotzero}
Let the TD system  $\Phi$ be as in 
(\ref{eq:phidisplay}).
Then for $0 \leq i \leq d$,
the space $U_i$ from Definition
\ref{def:meaningofVi}, the projection $F_i$ from
Definition
\ref{def:deffi}, and the scalar
$\rho_i$
from
Corollary \ref{thm:symunimodal} are all nonzero. 
\end{lemma}

\begin{proof} Recall $\rho_i$ is the dimension of $E_iV$, and this space is not zero, so $\rho_i\not=0$.
$U_i $ has dimension $\rho_i$, so $U_i\not=0$.
Recall $F_iV=U_i$, so $F_i\not=0$.

\end{proof}

\section{\bf The raising and lowering maps}

\noindent In this section, we continue to study the situation of
Definition
\ref{def:meaningofVi}. We introduce the raising and lowering
maps, and use these to show the scalars $\rho_i$ from
Corollary
\ref{thm:symunimodal} 
form a unimodal sequence. 

\begin{definition}
\label{def:defRandL}
With reference to 
Definition
\ref{def:meaningofVi},
we define
\begin{eqnarray}
R&=&A-\sum_{h=0}^d \theta_hF_h,
\label{eq:defR}
\\
L&=&A^*-\sum_{h=0}^d \theta^*_hF_h,
\label{eq:defL}
\end{eqnarray}
where  $F_0, F_1, \ldots, F_d$ are from
Definition
\ref{def:deffi}.
We refer to $R$ (resp. $L$) as the raising map (resp. lowering map). 
\end{definition}

\begin{lemma}
\label{lem:defRandLaction}
With reference to 
Definition
\ref{def:meaningofVi} and
Definition
\ref{def:defRandL}, for all integers $i$ $(0 \leq i \leq d)$, and for
all $v \in U_i$, 
\begin{equation}
Rv = (A-\theta_iI)v, \qquad \qquad Lv = (A^*-\theta^*_iI)v.
\label{eq:randLvsAAs}
\end{equation}
\end{lemma}

\begin{proof} To get the equation on the left,
apply both sides of
(\ref{eq:defR}) to
$v$, and evaluate the result using
(\ref{eq:fiui}), 
(\ref{eq:fiuj}). 
The equation on the right is similarily obtained.

\end{proof}

\begin{corollary}
\label{lem:RandLbasic}
With reference to 
Definition
\ref{def:meaningofVi} and
Definition
\ref{def:defRandL},
\begin{enumerate}
\item $RU_i \subseteq U_{i+1} \qquad \qquad (0 \leq i <d)$, $\qquad RU_d=0$. 
\item $LU_i \subseteq U_{i-1} \qquad \qquad  (0 < i \leq d)$, $\qquad LU_0=0$.
\end{enumerate}
\end{corollary}

\begin{proof} Combine Lemma
\ref{lem:defRandLaction} with 
(\ref{eq:Aaction}) and  
(\ref{eq:Asaction}).  

\end{proof}

\begin{lemma}
\label{lem:RandLnext}
With reference to Definition
\ref{def:meaningofVi},
Definition \ref{def:deffi},
and
Definition
\ref{def:defRandL},
\begin{enumerate}
\item $RF_i=F_{i+1}R \qquad \qquad (0 \leq i < d), \qquad RF_d = 0, \qquad F_0R = 0$. 
\item $LF_i=F_{i-1}L \qquad \qquad (0< i \leq d), \qquad  LF_0=0, \qquad F_dL = 0$.
\item $R^{d+1} = 0$.
\item $L^{d+1} = 0$.
\end{enumerate}
\end{lemma}

\begin{proof} (i) Pick any integer $i$ $(0 \leq i <d)$.
Using 
Definition
\ref{def:deffi}
and 
Corollary
\ref{lem:RandLbasic}(i), we find the expression
$RF_i - F_{i+1}R $ vanishes on each of $U_0, U_1,\ldots, U_d$,
so this expression is zero
in view of 
(\ref{eq:vdecomp}). Similarily, $RF_d $ and $F_0R$ vanish on each of 
 $U_0, U_1,\ldots, U_d$, so they are zero.

\noindent (ii) Similar to the proof of (i).

\noindent (iii) Immediate from 
Corollary \ref{lem:RandLbasic}(i) and 
(\ref{eq:vdecomp}).

\noindent (iv) Immediate from  
Corollary \ref{lem:RandLbasic}(ii) and 
(\ref{eq:vdecomp}).

\end{proof}

\begin{lemma}
\label{lem:Risnotzero}
With reference to Definition
\ref{def:meaningofVi} and
Definition
\ref{def:defRandL},
%
pick any integers $i,j$ $(0 \leq i\leq j\leq d)$. Then
the linear transformation 
\begin{equation}
{{U_i \quad \; \rightarrow \;\quad U_j}\atop {v \quad \rightarrow \quad R^{j-i}v}}
\label{eq:Rsur}
\end{equation}
is an injection if $i+j\leq d$, a bijection if $i+j=d$,
and a surjection if $i+j\geq d$.
The linear transformation 
\begin{equation}
{{U_j \;\quad \rightarrow \;\quad U_i}\atop {v \quad \rightarrow \quad L^{j-i}v}}
\label{eq:Lsur}
\end{equation}
is an injection if $i+j\geq d$, a bijection if $i+j=d$,
and  a surjection if $i+j\leq d$.

\noindent
(Caution: the maps 
(\ref{eq:Rsur}), 
(\ref{eq:Lsur})
are not inverses in general, even in the case  $i+j=d$). 
\end{lemma}

\begin{proof} Concerning the map
(\ref{eq:Rsur}), 
first assume $i+j \leq d$.
To show 
(\ref{eq:Rsur}) is an injection, 
 we pick any
vector $v \in U_i$
such that $R^{j-i}v=0$, and show $v=0$.
Using Lemma 
\ref{lem:defRandLaction}  and Corollary 
\ref{lem:RandLbasic}(i),
\beast
0 &=& R^{j-i}v
\\
&=& (A-\theta_{i}I)(A-\theta_{i+1}I)\cdots (A-\theta_{j-1}I)v,
\eeast
so
\begin{equation}
v \in E_iV+E_{i+1}V+\cdots +E_{j-1}V.
\label{eq:vinsum}
\end{equation}
By (\ref{eq:vsumzeroi}), and since
$v \in U_i$,    
\begin{equation}
v \in E^*_0V+E^*_1V+\cdots + E^*_iV.
\label{eq:vinFiv}
\end{equation}
Applying Lemma 
\ref{lem:manyvijzeroS99} to
 $\Phi^{\Downarrow}$, and using $j\leq d-i$, we find
\begin{equation}
(E^*_0V+E^*_1V+\cdots + E^*_iV)\cap (E_0V+E_1V+\cdots + E_{j-1}V)=0.
\label{eq:lem37again}
\end{equation}
Combining 
(\ref{eq:vinsum}),
(\ref{eq:vinFiv}),
(\ref{eq:lem37again}), we find $v=0$, and
it follows
(\ref{eq:Rsur}) is an injection.
Next suppose $i+j=d$. In this case $U_i$ and $U_j$ have the same
dimension by Corollary
\ref{thm:symunimodal},
so 
the injection (\ref{eq:Rsur}) is actually a bijection. 
We now assume $i+j\geq d$, and  show
(\ref{eq:Rsur}) is a surjection. To do this, 
we pick any
$w \in U_j$, 
and obtain an element $v \in U_i$ such that $R^{j-i}v = w$. 
From our comments above, the restriction of
$R^{2j-d}$ to $U_{d-j}$ is a bijection onto $U_j$, so
there exists $u \in U_{d-j}$ such that $R^{2j-d}u = w$.
Set $v = R^{i+j-d}u$. Applying Corollary 
\ref{lem:RandLbasic}(i), we find
$v \in U_i$. Also
\beast
R^{j-i}v  &=& R^{2j-d}u
\\
&=& w,
\eeast
and it follows (\ref{eq:Rsur}) is a surjection.  
We have now proved all our claims about 
(\ref{eq:Rsur}). To  obtain our claims about
(\ref{eq:Lsur}), apply the first part of 
the lemma to 
 $\Phi^{\downarrow \Downarrow *}$.

\end{proof}

\begin{corollary}
\label{cor:symunimodal} 
The scalars 
$\rho_0, \rho_1,\ldots  \rho_d$ 
from 
Corollary \ref{thm:symunimodal} 
satisfy 
$\rho_{i-1} \leq \rho_i$ for $1 \leq i \leq d/2$.  

\end{corollary}

\begin{proof}
Let the spaces $U_0, U_1,\ldots, U_d$ be as in 
Definition
\ref{def:meaningofVi} and recall $\rho_i$ is the dimension
of $U_i$ for $0 \leq i\leq d$. Pick any integer $i$ $(1\leq i\leq d/2)$.
By Lemma
\ref{lem:Risnotzero}, we find the restriction of the raising map
$R$ to
$U_{i-1} $ is an injection into $U_i$, so $\rho_{i-1}\leq \rho_i$.

\end{proof}

\noindent We end this section with a result we will use later in the paper.

\begin{corollary}
\label{cor:Rnotzero}
Let the TD system $\Phi$ be as in 
Definition
\ref{def:meaningofVi},
let the projections $F_0, F_1, \ldots, F_d$ be
as in 
Definition
\ref{def:deffi},
and let the maps $R,L$ be as in 
Definition
\ref{def:defRandL}.
Then for $0 \leq i \leq j\leq d$, we have $R^{j-i}F_i\not=0$
and $L^{j-i}F_j\not=0$.
\end{corollary}
\begin{proof} These maps are nonzero
by
Lemma 
\ref{lem:thevifinotzero}
and Lemma
\ref{lem:Risnotzero}.

\end{proof}

\section{Writing products involving $A$ and $A^*$ in terms of $R$ and $L$} 

\medskip
\noindent When working with a TD system such as
(\ref{eq:phidisplay}),  
instead of working
directly with the elements $A$ and $A^*$, it is often easier to work
with the correspondng raising map $R$ and  lowering map $L$ from
Definition
\ref{def:defRandL}.
In this section, we consider how to write products 
involving $A$ and $A^*$ in terms of $R$ and $L$. 
We begin with some simple observations.

\begin{lemma}
\label{lem:aandastarRL}
With reference to 
Definition
\ref{def:meaningofVi} and
Definition
\ref{def:deffi},
\begin{enumerate}
\item $F_iAF_i = \theta_i F_i\qquad \qquad (0 \leq i \leq d)$, 
\item $F_{i+1}AF_i = RF_i\qquad \qquad (0 \leq i < d)$, 
\item $F_jAF_i = 0 \quad \hbox{if}\quad  j-i \not\in  \lbrace 0,1\rbrace,
\qquad  \qquad (0 \leq i ,j\leq d)$.
\end{enumerate}
\end{lemma}

\begin{proof} To verify these equations, in each case
eliminate $A$ using
(\ref{eq:defR}), and evaluate the result using Lemma
\ref{lem:fielem} and Lemma
\ref{lem:RandLnext}.

\end{proof}

\begin{lemma}
\label{lem:astarandRL}
With reference to 
Definition
\ref{def:meaningofVi} and
Definition
\ref{def:deffi},
\begin{enumerate}
\item $F_iA^*F_i = \theta^*_i F_i\qquad \qquad (0 \leq i \leq d)$, 
\item $F_{i-1}A^*F_i = LF_i\qquad \qquad (0 < i \leq  d)$, 
\item $F_jA^*F_i = 0 \quad \hbox{if}\quad  j-i \not\in  \lbrace 0,-1\rbrace,
\qquad  \qquad (0 \leq i ,j\leq d)$.
\end{enumerate}
\end{lemma}

\begin{proof} Similar to the  proof of 
Lemma 
\ref{lem:aandastarRL}.

\end{proof}

\begin{lemma}
\label{lem:AandAstartoRL}
With reference to 
Definition
\ref{def:meaningofVi} and
Definition
\ref{def:deffi},
 let $n$ denote a positive
integer, and  
let $B_1, B_2, \ldots, B_n $
denote a sequence of
elements taken from the set $\lbrace A,A^*\rbrace $.
Then for $0 \leq r,s \leq d$, the  product 
\begin{equation}
F_rB_1B_2\cdots B_nF_s = 
\sum_\sigma
W_\sigma F_s
\label{eq:onewaytowritesum}
\end{equation}
with the sum interpreted as follows.
The sum 
 is over all sequences  of integers 
\begin{equation}
\sigma = (i_0, i_1, \ldots, i_n) 
\label{eq:sequenceinsum}
\end{equation}
such that 
$0 \leq i_j \leq  d$ for $0 \leq j \leq n$, and such that
\beast
&&\qquad i_0=r, \qquad i_n=s, 
\\
&&i_j - i_{j-1} \in  \lbrack B_j \rbrack \qquad \qquad (1 \leq  j \leq n), 
\eeast
where we define $\lbrack A \rbrack = \lbrace 0,-1 \rbrace $
and 
$\lbrack A^* \rbrack = \lbrace 0,1 \rbrace $.
For each sequence $\sigma $ in 
(\ref{eq:sequenceinsum}),
the corresponding $W_\sigma $ is the product 
$W_\sigma = W_1 W_2\cdots W_n$, where the $W_i$ are given by
\beast
W_j=\cases{ 
R, &if $\;    
i_j = i_{j-1} - 1$;\cr
L, &if $\;    
i_j = i_{j-1} + 1$;\cr
\theta_{i_j}, & if $\;
i_j = i_{j-1}$ and $B_j = A$; \cr
\theta^*_{i_j}, & if $\;
i_j = i_{j-1} $ and $ B_j = A^*$\cr} \qquad \qquad (1 \leq j \leq n),
\eeast
and where $R$ and $L$ are from
Definition
\ref{def:defRandL}.
\end{lemma}

\begin{proof} 
The left side of 
(\ref{eq:onewaytowritesum}) equals
$F_rB_1IB_2I\cdots IB_nF_s$.
Eliminating each copy of $I$ in this product 
using (\ref{eq:fisumtoident}), and evaluating the result using
Lemma
\ref{lem:aandastarRL} and Lemma
\ref{lem:astarandRL}, we obtain the sum on the
right side of 
(\ref{eq:onewaytowritesum}).

\end{proof}

\begin{example}
\label{ex:howtocompute}
For $0 \leq i \leq d-1$, we have
\beast
F_{i+1}AA^*AF_i = RLRF_i + \theta_{i+1}\theta^*_{i+1}RF_i
+ \theta_i\theta^*_iRF_i.
\eeast
\end{example}
\begin{proof} Routine application of Lemma
\ref{lem:AandAstartoRL}.

\end{proof}

\section{Recurrent sequences}

\noindent  It is going to turn out that the eigenvalue
sequence and dual eigenvalue sequence of a TD system
each satisfy a certain recurrence.
In this section, we set the stage by considering
this recurrence from several points of view.

\begin{definition}
\label{def:recseqsettupS99}
In this section, $\fld $ will denote a field,
$d$ will denote a nonnegative integer,  and 
 $\; \theta_0, \theta_1, \ldots,  \theta_d\; $  will denote a sequence
 of scalars taken from $\fld$.

\end{definition}

\begin{definition}
\label{lem:beginthreetermS99}
With reference to Definition
\ref{def:recseqsettupS99}, let $\beta, \gamma, \varrho$ denote scalars
in $\fld$.
\begin{enumerate} 
\item The sequence 
  $\; \theta_0, \theta_1, \ldots,  \theta_d \;$ 
is said to be recurrent whenever $\theta_{i-1}\not=\theta_i$ for
$2 \leq i \leq d-1$, and 
\begin{equation}
{{\theta_{i-2}-\theta_{i+1}}\over {\theta_{i-1}-\theta_i}} 
\label{eq:thethingwhichisbetaS99}
\end{equation}
is independent of
$i$, for $2 \leq i \leq  d-1$.
\item The sequence 
  $\; \theta_0, \theta_1, \ldots,  \theta_d \;$ 
is said to be $\beta$-recurrent whenever 
\begin{equation}
\theta_{i-2}\,-\,(\beta+1)\theta_{i-1}\,+\,(\beta +1)\theta_i \,-\,\theta_{i+1}
\label{eq:betarecS99}
\end{equation}
is zero for 
$2 \leq i \leq d-1$.
\item The sequence 
  $\; \theta_0, \theta_1, \ldots,  \theta_d \;$ 
is said to be $(\beta,\gamma)$-recurrent whenever 
\begin{equation}
\theta_{i-1}\,-\,\beta \theta_i\,+\,\theta_{i+1}=\gamma 
\label{eq:gammathreetermS99}
\end{equation}
 for 
$1 \leq i \leq d-1$.
\item The sequence 
  $\; \theta_0, \theta_1, \ldots,  \theta_d \;$ 
is said to be $(\beta,\gamma,\varrho)$-recurrent whenever 
\begin{equation}
\theta^2_{i-1}-\beta \theta_{i-1}\theta_i+\theta^2_i 
-\gamma (\theta_{i-1} +\theta_i)=\varrho
\label{eq:varrhothreetermS99}
\end{equation}
 for 
$1 \leq i \leq d$.
\end{enumerate}
\end{definition}

\begin{lemma} 
\label{lem:recvsbrecS99}
With reference to Definition
\ref{def:recseqsettupS99}, the following are equivalent.
\begin{enumerate}
\item The sequence 
  $\; \theta_0, \theta_1, \ldots,  \theta_d \;$ 
is  recurrent.
\item There exists $\beta \in \fld$ such that
  $\; \theta_0, \theta_1, \ldots,  \theta_d \;$  is $\beta$-recurrent,
and $\theta_{i-1}\not=\theta_i$ for
$2 \leq i \leq d-1$.
\end{enumerate}
\noindent Suppose (i), (ii), and that $d\geq 3$. Then
the common value of
(\ref{eq:thethingwhichisbetaS99}) 
equals $\beta +1$.

\end{lemma}
\begin{proof} Routine.

\end{proof}

\begin{lemma}
\label{lem:brecvsbgrecS99}
With reference to Definition
\ref{def:recseqsettupS99}, the
following are equivalent for all $\beta \in \fld$.
\begin{enumerate}
\item The sequence 
  $\; \theta_0, \theta_1, \ldots,  \theta_d \;$ 
is  $\beta$-recurrent.
\item There exists $\gamma \in \fld$ such that
  $\; \theta_0, \theta_1, \ldots,  \theta_d \;$  
  is $(\beta,\gamma)$-recurrent.
\end{enumerate}
\end{lemma}
\begin{proof} 
$(i)\rightarrow (ii) $ 
For $2\leq i \leq d-1$, the expression
(\ref{eq:betarecS99}) is zero by assumption,
so
\beast
\theta_{i-2}\,-\,\beta \theta_{i-1}\,+\,\theta_i \;= \;
\theta_{i-1}\,-\,\beta \theta_i\,+\,\theta_{i+1}.
\eeast
Apparently the left side of 
(\ref{eq:gammathreetermS99}) is independent of $i$, and
the result follows.

\noindent 
$(ii)\rightarrow (i) $  Subtracting the equation 
(\ref{eq:gammathreetermS99}) at $i$ from the corresponding equation
obtained by replacing $i$ by $i-1$, we find
(\ref{eq:betarecS99}) is zero 
for $2\leq i \leq d-1$.

\end{proof}

\begin{lemma}
\label{lem:bgrecvsbgdrecS99}
With reference to Definition
\ref{def:recseqsettupS99}, the following (i),(ii) hold
for all $\beta, \gamma \in \fld$.
\begin{enumerate}
\item  Suppose 
  $\; \theta_0, \theta_1, \ldots,  \theta_d \;$ 
is  $(\beta,\gamma)$-recurrent. Then 
there exists $\varrho \in \fld$ such that
  $\; \theta_0, \theta_1, \ldots,  \theta_d \;$  is $(\beta,\gamma,\varrho)$-recurrent.
\item  Suppose 
  $\; \theta_0, \theta_1, \ldots,  \theta_d \;$ 
is  $(\beta,\gamma,\varrho)$-recurrent, and that $\theta_{i-1}\not=\theta_{i+1}$
for $1 \leq i\leq d-1$. Then
  $\; \theta_0, \theta_1, \ldots,  \theta_d \;$  is $(\beta,\gamma)$-recurrent.
\end{enumerate}
\end{lemma}

\begin{proof} 
Let  $p_i$ denote the expression on the left in
(\ref{eq:varrhothreetermS99}),
and observe
 \beast
p_i-p_{i+1} &=& 
(\theta_{i-1}-\theta_{i+1})(\theta_{i-1}-\beta \theta_i +\theta_{i+1} - \gamma)
\eeast
for $1 \leq i \leq d-1$. 
Assertions (i), (ii) are both routine consequences of this.

\end{proof}

\noindent We mention the standard parametric expressions
for recurrent sequences.
To state the result, we 
let 
 ${\cal F}^{cl}$ denote the algebraic closure of $\fld$.
%
%
%
%
%
%
\begin{lemma}
\label{lem:closedformthreetermS99}
With reference to 
Definition
\ref{def:recseqsettupS99},
 pick any $\beta
 \in \fld$, and assume  
 $\; \theta_0, \theta_1, \ldots,  \theta_d\; $ is $\beta$-recurrent. 
Then
the following  (i)--(iv) hold.
\begin{enumerate}
\item  Suppose $\beta \not=2$, $\beta \not=-2$, and pick
$q \in 
{\cal F}^{cl}$ such that 
 $q+q^{-1}=\beta $. Then there exists  scalars 
 $\alpha_1, \alpha_2, \alpha_3$  in  
${\cal F}^{cl}$ such that 
\begin{equation}
\theta_i = \alpha_1 + \alpha_2 q^i + \alpha_3 q^{-i}
\qquad \qquad (0 \leq i \leq d). 
\label{eq:closedformthreetermIS99}
\end{equation}
\item Suppose $\beta = 2$ and $\hbox{char}(\fld) \not=2$. Then there
exists 
 $\alpha_1, \alpha_2, \alpha_3 $  in $\fld $ such that
\begin{equation}
\theta_i = \alpha_1 + \alpha_2 i + \alpha_3 i^2 
\qquad \qquad (0 \leq i \leq d).  
\label{eq:closedformthreetermIIS99}
\end{equation}
\item Suppose $\beta = -2$ and  $\hbox{char}(\fld) \not=2$. Then there
exists 
 $\alpha_1, \alpha_2, \alpha_3 $  in $\fld $ such that
\begin{equation}
\theta_i = \alpha_1 + \alpha_2 (-1)^i + \alpha_3 i(-1)^i 
\qquad \qquad (0 \leq i \leq d).  
\label{eq:closedformthreetermIIIS99}
\end{equation}
\item Suppose $\beta = 0$ and  $\hbox{char}(\fld) =2$. Then there
exists 
 $\alpha_1, \alpha_2, \alpha_3 $  in $\fld $ such that
\begin{equation}
\theta_i = \alpha_1 + \alpha_2 i + \alpha_3 \Biggl( {{i }\atop {2}}\Biggr) 
\qquad \qquad (0 \leq i \leq d),  
\label{eq:closedformthreetermIVS99}
\end{equation}
where we define 
\beast
 \Biggl( {{i }\atop {2}}\Biggr) 
 = \left\{ \begin{array}{ll}
                   0  & \mbox{if $\;i=0\;$ or $\;i=1\;$ (mod $4$), } \\
			  1 & \mbox{if $\;i=2 \;$ or $\;i=3\;$ (mod $4$). }
			   \end{array}
			\right. 
\eeast
\end{enumerate}
\end{lemma}

%
%
%

\section{A few comments}

\medskip
\noindent In the next section,
we are going to show the elements $A$ and $A^*$ in a TD pair
satisfy two cubic polynomial equations. In this section, we obtain
a few preliminary results that we will need.

\begin{definition}
\label{def:prelimsit}
Throughout this section,
 $\fld$ will denote a field, and $V$ will denote a vector space
over $\fld$ with finite positive dimension. We 
let $A$ and $A^*$
denote linear transformations on $V$, with $A$ 
diagonalizable. We let $E_0, E_1, \ldots, E_d$ denote
an ordering of the primitive idempotents of $A$.
We let $\theta_i$ denote the
eigenvalue of $A$ associated with $E_i$, for
$0 \leq i \leq d$.
We remark $\theta_0, \theta_1, \ldots, \theta_d$ are mutually distinct.

\end{definition}

\begin{lemma}
\label{lem:caldpre}
With reference to 
Definition
\ref{def:prelimsit},
suppose 
\begin{equation}
E_iA^*E_j = 0 \quad \hbox{if} \quad |i-j|>1, \qquad \qquad 
(0 \leq i,j\leq   d),
\label{eq:caldpre1}
\end{equation}
and let 
${\cal D}$ denote the subalgebra of $\hbox{End}(V)$ generated by $A$.
Then
\begin{equation}
\Span{XA^*Y-YA^*X \;|\;X, Y \in {\cal D}} \,=\,\lbrace XA^*-A^*X\;|\;X \in 
{\cal D}\rbrace. \qquad \quad 
\label{eq:calDS99n}
\end{equation}

\end{lemma}

\begin{proof}
 For notational convenience set $E_{-1}=0$, $E_{d+1}=0$.
Observe by
(\ref{eq:primid3S99})
and
(\ref{eq:caldpre1}) 
that for  $0 \leq j \leq d$, 
both
\begin{eqnarray}
A^*E_j&=& E_{j-1}A^*E_j+E_jA^*E_j + E_{j+1}A^*E_j,
\label{eq:Aexpand2S99n}
\\
E_jA^*&=&E_{j}A^*E_{j-1}+E_jA^*E_j + E_{j}A^*E_{j+1}.
\label{eq:Aexpand1S99n}
\end{eqnarray}
Pick any integer $i$ $(0 \leq i \leq d)$.
Summing  
(\ref{eq:Aexpand2S99n})
over 
 $j=0,1,\ldots,i$, summing 
(\ref{eq:Aexpand1S99n}) over 
 $j=0,1,\ldots,i$, and taking the difference between the
 two sums, 
we obtain 
\begin{equation}
E_iA^*E_{i+1} - E_{i+1}A^*E_i \;=\;L_i A^* - A^* L_i,
\label{eq:liconversionS99n}
\end{equation}
where  $L_i := E_0+E_1+\cdots +E_i $.
Observe 
$\cal D$ is spanned by both $E_0, E_1, \ldots, E_d$
and $L_0, L_1, \ldots,L_d$, so 
\beast
&&\Span{XA^*Y-YA^*X \;|\;X, Y \in {\cal D}}
\\
&& \qquad \qquad =\; \Span{E_iA^*E_j-E_jA^*E_i\;|\;
0 \leq i,j\leq d}\\
 && \qquad \qquad   
=\; \Span{E_iA^*E_{i+1}-E_{i+1}A^*E_i\;|\;
0 \leq i\leq d}\\
 &&\qquad \qquad =\; \Span{L_iA^*-A^*L_i\;|\;0 \leq i \leq d} \\
&&\qquad \qquad =\;\lbrace XA^*-A^*X\;|\;X \in 
{\cal D}\rbrace,
\eeast
and we have
(\ref{eq:calDS99n}).

\end{proof}

\begin{lemma}
\label{lem:caldpre2}
With reference to 
Definition
\ref{def:prelimsit},
suppose 
\begin{equation}
E_iA^*E_j = \cases{0, &if $\;\vert i-j\vert > 1$;\cr
\not=0, &if $\;\vert i-j \vert = 1$\cr}
\qquad \qquad 
(0 \leq i,j\leq d),
\label{eq:caldpre2}
\end{equation}
and pick any scalars $\beta, \gamma, \varrho $ in $\fld$.
Then
\begin{eqnarray}
0 &=&\lbrack A,A^2A^*-\beta AA^*A + 
A^*A^2 -\gamma (AA^*+A^*A)-\varrho A^*\rbrack 
\label{eq:qdolanA}
\end{eqnarray}
if and only if the sequence $\theta_0, \theta_1, \ldots, \theta_d$
is $(\beta, \gamma, \varrho)$-recurrent.

\end{lemma}

\begin{proof} For notational convenience,  
we define a two variable polynomial
$p \in \fld\lbrack \lambda,\mu \rbrack $
by
\begin{equation}
p(\lambda, \mu) = \lambda^2 - \beta \lambda \mu + \mu^2 - \gamma (\lambda + \mu)- \varrho.
\label{eq:defoftwovarp}
\end{equation}
Let $C$ denote the expression on the right in
(\ref{eq:qdolanA}), and  observe 
\begin{equation}
C = \sum_{i=0}^d \sum_{j=0}^d E_iCE_j
\label{eq:csum}
\end{equation}
in view of 
(\ref{eq:primid3S99}).
Using 
(\ref{eq:primid1S99}) and the definition of $C$, we find that for $0 \leq i,j\leq d$,
\begin{equation}
E_iCE_j=E_iA^*E_j
(\theta_i-\theta_j)p(\theta_i,\theta_j).
\label{eq:eicej}
\end{equation}
Suppose
(\ref{eq:qdolanA}) holds, so that    
$C=0$. Then for $1\leq i \leq d$, we have
 $E_{i-1}CE_i=0$,
and also $E_{i-1}A^*E_i\not=0$ by
(\ref{eq:caldpre2}), 
so $p(\theta_{i-1},\theta_i)=0$ by
(\ref{eq:eicej}).
Evaluating this using
(\ref{eq:defoftwovarp})  and Definition
\ref{lem:beginthreetermS99}(iv), 
we find
$\theta_0, \theta_1, \ldots, \theta_d$ is 
$(\beta, \gamma, \varrho)$-recurrent.
Conversely, suppose 
$\theta_0, \theta_1, \ldots, \theta_d$ is 
$(\beta, \gamma, \varrho)$-recurrent.
By
  Definition
\ref{lem:beginthreetermS99}(iv) and 
(\ref{eq:defoftwovarp}),
 we find  
  $p(\theta_{i-1}, \theta_i)
  =p(\theta_{i}, \theta_{i-1}) = 0$
for $1 \leq i \leq d$.  Now for $0 \leq i,j\leq d$,
at least one of the factors on the right in 
(\ref{eq:eicej}) equals zero, so 
$E_iCE_j=0$.
Apparently each term on the right in 
(\ref{eq:csum}) is zero, so $C=0$, and 
we have (\ref{eq:qdolanA}).

\end{proof}

\section{Two polynomial equations for $A$ and $A^*$}

\begin{theorem}
\label{eq:lastchancedolangradyS99}
Let $\fld$ denote a field, and let
 $(A,A^*)$ denote a TD pair over  $\fld$.
 Then
there  exists a sequence of scalars $\beta, \gamma, \gamma^*,
\varrho, \varrho^*$ taken from $\fld$ such that
both
\begin{eqnarray}
0 &=&\lbrack A,A^2A^*-\beta AA^*A + 
A^*A^2 -\gamma (AA^*+A^*A)-\varrho A^*\rbrack, 
\label{eq:qdolangrady199n}
\\
0 &=& \lbrack A^*,A^{*2}A-\beta A^*AA^* + AA^{*2} -\gamma^* (A^*A+AA^*)-
\varrho^* A\rbrack.  \qquad  \quad 
\label{eq:qdolangrady2S99n}
\end{eqnarray}
The sequence is unique if the diameter of the pair is at least 3.
\end{theorem}

\begin{proof} Let $\ls$ denote a TD system associated with $(A,A^*)$.
Let 
$\theta_0,\theta_1,\ldots, \theta_d$
denote the  
eigenvalue sequence of $\Phi$, and let 
$\theta^*_0,\theta^*_1,\ldots, \theta^*_d$ denote the 
dual eigenvalue sequence of $\Phi$.
Let $F_0, F_1, \ldots, F_d$ denote the projection
maps for $\Phi$, 
from Definition
\ref{def:deffi}.
Let $R$ and $L$ denote the raising  and lowering maps for $\Phi$,
from Definition
\ref{def:defRandL}.
First assume $d\geq 3$.
Observe $\Phi$ satisfies
(\ref{eq:caldpre1})
by
Definition
\ref{def:deflstalkS99}(iv), so Lemma 
\ref{lem:caldpre} applies.
Taking  $X=A^2, Y=A$ on the left in (\ref{eq:calDS99n}), we find
there exists scalars $\alpha_1, \alpha_2,\ldots,\alpha_d$ in $\fld$ such that
\begin{equation}
A^2A^*A - A A^*A^2 \;= \; \sum_{i=1}^d \alpha_i(A^iA^*-A^*A^i).
\label{eq:qaskeywilsonCS99n}
\end{equation}
We show $\alpha_i=0\;$ for $4 \leq i \leq d$. 
Suppose not, and set
\beast
t:= \hbox{max}\lbrace i \;|\;4 \leq i \leq d,\;\;\alpha_i \not=0\rbrace.
\eeast
Multiplying each term in 
(\ref{eq:qaskeywilsonCS99n}) on the left by $F_t$ and on the right
by $F_0$,  and evaluating the result using
Lemma \ref{lem:AandAstartoRL},
we routinely find $R^tF_0$ times 
 $(\theta^*_0-\theta^*_t)\alpha_t$ is zero.
 However   
$R^tF_0\not=0$ by 
Corollary \ref{cor:Rnotzero},
$\theta^*_0\not=\theta^*_t$ by construction, and 
we assumed 
$\alpha_t\not=0$, 
so 
 we  have a contradiction.
We now have  
 $\,\alpha_i=0\, $ for $4 \leq i \leq d$,
 so (\ref{eq:qaskeywilsonCS99n}) becomes
\begin{eqnarray}
\displaystyle{
{{A^2A^*A-AA^*A^2
= \alpha_1(AA^*-A^*A)\;+\;\alpha_2(A^2A^*-A^*A^2)
}\atop
{
\qquad \qquad \qquad \qquad \qquad +\;
\alpha_3(A^3A^*-A^*A^3).}}}
\label{eq:qaskeywilsonDS99n}
\end{eqnarray}
We show  $\alpha_3\not=0$. Suppose $\alpha_3=0$. 
Multiplying each term in 
(\ref{eq:qaskeywilsonDS99n}) on the 
left  
by $F_3$
and on the right by $F_0$,
and evaluating the result using Lemma 
\ref{lem:AandAstartoRL},
we routinely find $R^3F_0$ times
$\theta^*_1-\theta^*_2$
is zero.  Observe 
$R^3F_0\not=0$  by
Corollary \ref{cor:Rnotzero}, 
and  
$\theta^*_1\not=\theta^*_{2}$,
so we have a contradiction. We have now shown $\alpha_3\not=0$.
Set 
\begin{eqnarray}
\displaystyle{
{{C:= \alpha_1 A^*\,+\,\alpha_2(AA^*+A^*A)\,+\,\alpha_3(A^2A^*+A^*A^2)
\qquad \qquad \qquad 
}\atop {\qquad  \qquad \qquad \qquad \qquad \qquad \qquad  \,+\;
(\alpha_3-1)AA^*A.}}}
\label{eq:qaskeywilsonFS99n}
\end{eqnarray}
Observe $AC-CA$ equals  
\beast
&&\alpha_1(AA^*-A^*A)\;+\;\alpha_2(A^2A^*-A^*A^2)\;+\;
\alpha_3(A^3A^*-A^*A^3)
\\
&& \qquad \qquad \qquad \qquad \qquad \qquad \;+\;
AA^*A^2-A^2A^*A,
\eeast
and this equals zero  
in view of (\ref{eq:qaskeywilsonDS99n}).
Hence $A$ and $C$ commute.
Dividing 
$C$ by
 $\alpha_3$ and using
(\ref{eq:qaskeywilsonFS99n}), 
 we find $A$ commutes with
\beast
A^2A^*-\beta AA^*A + A^*A^2 -\gamma (AA^*+A^*A)-\varrho A^*,
\eeast
where
\begin{equation}
\beta := \alpha_3^{-1} - 1,\qquad \gamma := -\alpha_2 \alpha_3^{-1},
\qquad \varrho := -\alpha_1\alpha_3^{-1}.
\label{eq:defbetagam}
\end{equation}
We now have
(\ref{eq:qdolangrady199n}).
Pick any integer $i$ $(2 \leq i \leq d-1)$. 
Multiplying each term in
(\ref{eq:qdolangrady199n})
 on the right by $F_{i-2}$ and on the left  
by $F_{i+1}$, and evaluating the result using Lemma 
\ref{lem:AandAstartoRL},
we routinely find $R^3F_{i-2}$ times
\begin{equation}
\theta^*_{i-2}-(\beta+1)\theta^*_{i-1} + (\beta+1)\theta^*_{i} - \theta^*_{i+1}
\label{eq:alpha3recbet}
\end{equation}
is zero. 
 Observe 
 $R^3F_{i-2}\not=0$  by
Corollary \ref{cor:Rnotzero}, 
so 
(\ref{eq:alpha3recbet}) is zero, 
and it follows 
$\theta^*_0, \theta^*_1,\ldots, \theta^*_d$ is $\beta$-recurrent
in view of
Definition
\ref{lem:beginthreetermS99}(ii).
Applying
Lemma \ref{lem:brecvsbgrecS99} and then
Lemma \ref{lem:bgrecvsbgdrecS99} to the sequence
$\theta^*_0, \theta^*_1,\ldots, \theta^*_d$,  
we find there exists scalars $\gamma^*, \varrho^*$ in $\fld$
such that  this sequence is 
$(\beta,\gamma^*,\varrho^*)$-recurrent.
Applying Lemma
\ref{lem:caldpre2} to $\Phi^*$, we obtain
(\ref{eq:qdolangrady2S99n}).
Concerning uniqueness, 
let $\beta, \gamma, \gamma^*, \varrho, \varrho^*$ denote any scalars
in $\fld$ satisfying
(\ref{eq:qdolangrady199n}),
(\ref{eq:qdolangrady2S99n}).
Applying Lemma 
\ref{lem:caldpre2},
we find 
$\theta_0, \theta_1,\ldots, \theta_d$ is $(\beta,\gamma,\varrho)$-recurrent.
From this and our assumption $d\geq 3$,
one routinely finds 
 $\beta,\gamma,\varrho$ are uniquely determined.
Similarily $\gamma^*$ and $\varrho^*$ are uniquely determined, and
we have proved the theorem for the case  $d\geq 3$.
 Next assume $d\leq 2$,
and let $\beta $ denote any scalar in $\fld $.
If $d=2$ set
$\gamma = \theta_0-\beta \theta_1+\theta_2$,
and if $d\leq 1$ let $\gamma $
denote any scalar in $\fld$.
If $d\geq 1$ set 
\beast
\varrho=\theta^2_0-\beta \theta_0\theta_1+\theta_1^2-\gamma(\theta_0+\theta_1),
\eeast
and if $d=0$ let $\varrho$ denote any scalar in $\fld$.
Applying Definition
\ref{lem:beginthreetermS99}(iv), we find
the  sequence $\theta_0,\theta_1,\ldots, \theta_d$ is 
 $(\beta,\gamma,\varrho)$-recurrent.
  Applying 
Lemma 
\ref{lem:caldpre2},  we find
(\ref{eq:qdolangrady199n}) holds.
Applying the above argument to $\Phi^*$, we find there exists scalars
$\gamma^*, \varrho^*$ for which 
(\ref{eq:qdolangrady2S99n}) holds.
\end{proof}

\begin{remark}
\label{rem:mastervsDG}
The Dolan-Grady relations 
(\ref{eq:DG1}), 
(\ref{eq:DG2}) are the same thing as 
(\ref{eq:qdolangrady199n}),
(\ref{eq:qdolangrady2S99n}) with 
$\beta = 2$, $\gamma = \gamma^*=0$, $\varrho=b^2$, 
$\varrho^*=b^{*2}$.
The quantum Serre relations
(\ref{eq:Serre1}), 
(\ref{eq:Serre2})
are the same thing as  (\ref{eq:qdolangrady199n}),
(\ref{eq:qdolangrady2S99n}) with 
$\beta = q+q^{-1}$, $\gamma = \gamma^*=0$, $\varrho=\varrho^*=0$. 
\end{remark}

\section{The eigenvalues and dual eigenvalues of a TD system}

\noindent In the proof of 
Theorem
\ref{eq:lastchancedolangradyS99},
 we got an indication  that 
 the eigenvalues and dual eigenvalues of a TD system are recurrent.
In this section, we give detailed results  along this line. 

\begin{theorem}
\label{eq:lastchancedolangradyS992}
Let $\fld$ denote a field,  
let $(A,A^*)$ denote a TD pair over  
 $\fld$, 
and let 
 $\beta, \gamma, \gamma^*,
\varrho, \varrho^*$ denote scalars in $\fld$ that satisfy
(\ref{eq:qdolangrady199n}), 
(\ref{eq:qdolangrady2S99n}). Let 
$\theta_0, \theta_1, \ldots, \theta_d$ (resp.
 $\theta^*_0, \theta^*_1, \ldots, \theta^*_d$).
denote the eigenvalue sequence (resp.
dual eigenvalue sequence), for an associated TD system.
Then the expressions
\begin{equation}
{{\theta_{i-2}-\theta_{i+1}}\over {\theta_{i-1}-\theta_i}},\qquad \qquad  
 {{\theta^*_{i-2}-\theta^*_{i+1}}\over {\theta^*_{i-1}-\theta^*_i}} 
 \qquad  \qquad 
\label{eq:expsbetaplusone}
\end{equation} 
 both equal $\beta +1$, for $\;2\leq i \leq d-1$.  
Moreover,
\begin{eqnarray}
\gamma &=& \theta_{i-1}-\beta \theta_i + \theta_{i+1} \qquad \qquad (1 \leq i \leq d-1),
\label{eq:gamform}
\\
\gamma^* &=& \theta^*_{i-1}-\beta \theta^*_i + \theta^*_{i+1} \qquad \qquad (1 \leq i \leq d-1),
\label{eq:gamsform}
\\
\varrho&=& \theta^2_{i-1}-\beta \theta_{i-1}\theta_i+\theta_i^2-\gamma (\theta_{i-1}+\theta_i) \qquad \qquad (1 \leq i \leq d),
\label{eq:vrhoform}
\\
\varrho^*&=& \theta^{*2}_{i-1}-\beta \theta^*_{i-1}\theta^*_i+\theta_i^{*2}-
\gamma^* (\theta^*_{i-1}+\theta^*_i) \qquad \qquad (1 \leq i \leq d).
\label{eq:vrhosform}
\end{eqnarray}

\end{theorem}

\begin{proof} Let $\Phi$ denote the TD system  referred to in the statement of
the theorem.
Applying Lemma 
\ref{lem:caldpre2} to $\Phi$, we find
 the sequence 
 $\theta_0, \theta_1, \ldots, \theta_d$ is 
 $(\beta, \gamma, \varrho)$-recurrent, giving  
(\ref{eq:vrhoform}).
Applying Lemma
\ref{lem:bgrecvsbgdrecS99}, and since the $\theta_i$ 
 are distinct,
we find the sequence
 $\theta_0, \theta_1, \ldots, \theta_d$ is 
 $(\beta, \gamma)$-recurrent, giving   
(\ref{eq:gamform}). 
Applying 
Lemma
\ref{lem:brecvsbgrecS99}, we find
 $\theta_0, \theta_1, \ldots, \theta_d$ is 
 $\beta$-recurrent, so 
the expressions on the left in
(\ref{eq:expsbetaplusone})
all equal $\beta+1$.
Applying our above comments to $\Phi^*$, we obtain the 
remaining assertions.

\end{proof}

\begin{theorem}
\label{thm:paramsolsforeseq}
Let $\fld$ denote a field. Let 
 $\Phi$ denote a TD system over  $\fld$, with 
  eigenvalue
 sequence
 $\theta_0, \theta_1, \ldots, \theta_d$ and dual  
 eigenvalue
 sequence
 $\theta^*_0, \theta^*_1, \ldots, \theta^*_d$.
Let $\beta $ denote the scalar from
Theorem
\ref{eq:lastchancedolangradyS99}. Then
(i)--(iv) hold below.
\begin{enumerate}
\item  
Suppose $\beta \not=2, \beta \not=-2$, and pick 
$q \in 
{\cal F}^{cl}$ such that 
 $q+q^{-1}=\beta $. Then there exists  scalars 
 $\alpha_1, \alpha_2, \alpha_3,
 \alpha^*_1, \alpha^*_2, \alpha^*_3$
 in  
${\cal F}^{cl}$ such that
\begin{eqnarray}
\theta_i &=& \alpha_1 + \alpha_2q^i + \alpha_3q^{-i},
\\
\theta^*_i &=& \alpha^*_1 + \alpha^*_2q^i + \alpha^*_3q^{-i},
\end{eqnarray}
for $0 \leq i \leq d$.
Moreover $q^i\not=1$ for $1 \leq i \leq d$.
\item Suppose $\beta = 2$ and $\hbox{char}(\fld) \not=2$. Then there
exists 
 $\alpha_1, \alpha_2, \alpha_3, 
 \alpha^*_1, \alpha^*_2, \alpha^*_3 $  in $\fld $ such that
\begin{eqnarray}
\theta_i &=& \alpha_1 + \alpha_2 i + \alpha_3 i^2,
\\
\theta^*_i &=& \alpha^*_1 + \alpha^*_2 i + \alpha^*_3 i^2, 
\end{eqnarray}
for $0 \leq i \leq d$.  
Moreover $\hbox{char}(\fld)=0$
or $\hbox{char}(\fld)>d$.
\item Suppose $\beta = -2$ and  $\hbox{char}(\fld) \not=2$. Then there
exists 
 $\alpha_1, \alpha_2, \alpha_3, 
 \alpha^*_1, \alpha^*_2, \alpha^*_3 $  in $\fld $ such that
\begin{eqnarray}
\theta_i &=& \alpha_1 + \alpha_2 (-1)^i + \alpha_3 i(-1)^i, 
\\
\theta^*_i &=& \alpha^*_1 + \alpha^*_2 (-1)^i + \alpha^*_3 i(-1)^i, 
\end{eqnarray}
for $0 \leq i \leq d$.  
Moreover $\hbox{char}(\fld)=0$
or $\hbox{char}(\fld)>d/2$.
\item Suppose $\beta = 0$ and  $\hbox{char}(\fld) =2$. Then $d\leq 3$.
\end{enumerate}

\end{theorem}

\begin{proof} Most of the assertions are immediate from
Lemma 
\ref{lem:closedformthreetermS99}, but a few points need
explanation.
Suppose we are in the situation of (i). 
 Then 
using
(\ref{eq:closedformthreetermIS99}), we find that for $1 \leq i \leq  d$,
$q^i=1 $ implies $\theta_i=\theta_0$, a contradiction. 
Suppose we are in the situation of (ii). 
Then using (\ref{eq:closedformthreetermIIS99}), we find
that for $1\leq i \leq d$, $i=0$ (in $\fld$) implies 
 $\theta_i=\theta_0$, a contradiction.
Hence 
$\hbox{char}(\fld)=0$
or $\hbox{char}(\fld)>d$.
Suppose we are in the situation of (iii). 
Then using  
(\ref{eq:closedformthreetermIIIS99}), we find that for
  $1 \leq i \leq d$,    
$i=0$ (in $\fld $) and $i$ even implies 
 $\theta_i=\theta_0$.
Hence
 $\hbox{char}(\fld)=0$ or 
 $\hbox{char}(\fld)>d/2$. 
Suppose we are in the situation of (iv), and that  
 $d\geq 4$.
Then
applying  
(\ref{eq:closedformthreetermIVS99}), we find $\theta_0=\theta_4$, 
a contradiction.
Hence $d\leq 3$.
\end{proof}

\section{The raising and lowering maps, revisited} 

\noindent In this section we consider the implications of
Theorem 
\ref{eq:lastchancedolangradyS99}
for the raising map $R$ and the lowering map $L$.

\begin{theorem}
\label{thm:howRLinteract}
Let $\Phi$ denote a TD system,
with eigenvalue sequence
 $\theta_0, \theta_1, \ldots, \theta_d$ and dual eigenvalue sequence  
 $\theta^*_0, \theta^*_1, \ldots, \theta^*_d$.
Let $U_0, U_1,\ldots, U_d$ denote the spaces from 
Definition
\ref{def:meaningofVi}. Let  
$R$  and $L$ denote the raising and lowering maps for $\Phi$, 
from Definition
\ref{def:defRandL}.
Then for $0 \leq i \leq d-2$,
\begin{equation}
R^3L - (\beta+1)R^2LR+(\beta+1)RLR^2-LR^3 + (\beta+1)\varepsilon_i R^2
\label{eq:RLint1}
\end{equation}
vanishes on $U_i$, and
\begin{equation}
RL^3 - (\beta+1)LRL^2+(\beta+1)L^2RL-L^3R + (\beta+1)\varepsilon_i L^2
\label{eq:RLint2}
\end{equation}
vanishes on $U_{i+2}$, 
where $\beta $ is from
Theorem
\ref{eq:lastchancedolangradyS99}, and where
\begin{equation}
\varepsilon_i = (\theta_i-\theta_{i+2})(\theta^*_{i+1}-\theta^*_{i+2})-
(\theta^*_{i+2}-\theta^*_{i})(\theta_{i+1}-\theta_{i}).
\label{eq:defofvareps}
\end{equation}
\end{theorem}

\begin{proof}  
Let the projections $F_0, F_1, \ldots, F_d$ be as in 
Definition
\ref{def:deffi}.
To get 
(\ref{eq:RLint1}), multiply each term in 
(\ref{eq:qdolangrady199n}) on the right by $F_i$,
and on the left by $F_{i+2}$, and simplify using
Lemma
\ref{lem:AandAstartoRL}.
 To get 
(\ref{eq:RLint2}), multiply each term in 
(\ref{eq:qdolangrady2S99n}) on the left by $F_i$, and on the right by $F_{i+2}$, 
and simplify using
Lemma
\ref{lem:AandAstartoRL}.

\end{proof}

\noindent Referring to the above theorem, it is natural to consider when 
$\varepsilon_i$ is independent
of $i$. Below, we give one instance where this  occurs. 

\begin{theorem}
\label{thm:howRLinteractnice}
Let $\fld$ denote a field.
Let $\Phi$ denote a TD system over $\fld$,
with eigenvalue sequence
$\theta_0, \theta_1, \ldots, \theta_d$ and dual eigenvalue sequence  
 $\theta^*_0, \theta^*_1, \ldots, \theta^*_d$.
Suppose  
 there exists scalars 
$q, a,b,a^*, c^* $ in $\fld$  such that  
\begin{equation}
\theta_i = a+bq^i, \qquad \qquad \theta^*_i = a^*+ c^*q^{-i} \qquad \qquad (0 \leq i \leq d).
\label{eq:niceeigs}
\end{equation}
Then the raising map
$R$  and the lowering map $L$
from Definition
\ref{def:defRandL} satisfy the quantum Serre relations
\begin{equation}
0=\lbrack R,  R^2L - (q+q^{-1})RLR+LR^2\rbrack,
\label{eq:RLint1nice}
\end{equation}
\begin{equation}
0 = \lbrack L, L^2R - (q+q^{-1})LRL+RL^2\rbrack,
\label{eq:RLint2nice}
\end{equation}
where $\lbrack r,s\rbrack $ means $rs-sr$.
\end{theorem}

\begin{proof} Let $\Delta $ denote the expression on the right in 
(\ref{eq:RLint1nice}).
Evaluating  
(\ref{eq:defofvareps}) using  
(\ref{eq:niceeigs}), we find
$\varepsilon_i= 0 $ for $0 \leq i \leq d-2$.
Evaluating  
(\ref{eq:expsbetaplusone}) using 
(\ref{eq:niceeigs}), we find $\beta = q+q^{-1}$.
Using this information, we find
(\ref{eq:RLint1})  equals $\Delta $. 
Applying Theorem
\ref{thm:howRLinteract},
we find 
$\Delta$ vanishes on  
$U_i$ for $0 \leq i \leq d-2$.
Using 
Corollary \ref{lem:RandLbasic}, we 
find $\Delta$ vanishes on $U_i$ for $i = d-1$ and   $i=d$.
Apparently $\Delta$ vanishes on each of
$U_0, U_1, \ldots,U_d$, and these span $V$, so $\Delta=0$, and
we have (\ref{eq:RLint1nice}). 
Applying a similar argument,
we obtain
(\ref{eq:RLint2nice}).

\end{proof}

\section{Remarks and conjectures}
\medskip
\noindent  In this section, we give some 
open problems and ideas for future research.

\medskip
\noindent As we have indicated earlier,
the main problem that we would like to see
solved is the following. 

\begin{problem} Classify all the TD pairs.
\end{problem}

\medskip
\noindent The above problem might be difficult, 
but the fact that there is a classification in the thin case
\cite{LS99} gives us hope. 
To make headway on the non-thin case, we suggest  one of the following two
problems.

\begin{problem}
Classify all the TD pairs for which the sequence
$\rho_0, \rho_1, \ldots, \rho_d$ from
Corollary \ref{thm:symunimodal} 
is of the form
$1,2,2,2,\ldots, 2,1$.
We remark these TD pairs play an important role in the theory
of $P$-and $Q$-polynomial schemes\cite{HobIto},\cite{Tan}.
\end{problem}

\begin{problem}
Classify all the TD systems that satisfy the assumption
of Theorem
\ref{thm:howRLinteractnice}, where
we assume the scalar $q$ in that theorem is  not a root
of unity.
Given the results of that
theorem, it is natural to guess these TD sytems are
related to representations of the quantum
affine algebra
$U_r({\widehat {sl}}_2)$, where $r^2=q$.
See \cite{CPqaa}, \cite{CPqaar} for  information on this algebra.
\end{problem}

\medskip
\noindent In light of Theorem
\ref{eq:lastchancedolangradyS99},
we consider the following generalization of a TD pair.
Let  $\fld$  denote a field, and let
$V$ denote a vector space over $\fld$ with finite positive
dimension.
By a {\it generalized TD pair} on $V$, we mean an ordered
pair $(A,A^*)$, where  
 $A$ and $A^*$ are  linear transformations on $V$ 
such that (i), (ii) hold below.
\begin{enumerate}
\item There exists scalars
 $\beta, \gamma, \gamma^*,
\varrho, \varrho^*$ taken from $\fld$ such that
\begin{eqnarray}
0 &=&\lbrack A,A^2A^*-\beta AA^*A + 
A^*A^2 -\gamma (AA^*+A^*A)-\varrho A^*\rbrack, 
\label{eq:qdolangradyfin1}
\\
0 &=& \lbrack A^*,A^{*2}A-\beta A^*AA^* + AA^{*2} -\gamma^* (A^*A+AA^*)-
\varrho^* A\rbrack,  \qquad  \quad 
\label{eq:qdolangradyfin2}
\end{eqnarray}
where $\lbrack r,s\rbrack $ means $rs-sr$.
\item There is no subspace $W$ of $V$ such  that  both $AW\subseteq W$,
$A^*W\subseteq W$, other than $W=0$ and $W=V$.
\end{enumerate}
By Theorem
\ref{eq:lastchancedolangradyS99}, any TD pair is a generalized
TD pair. However,  
we should not expect every generalized TD pair to 
be a TD pair.
 We saw some indications of this 
in Example
\ref{ex:sl2easy},  
Example \ref{def:onsager}, 
and Example \ref{def:serre}.
In  each of these examples, to
get  a TD pair,
we needed to make an 
 assumption about the ground field $\fld$.
Also,  in Example
 \ref{def:serre}, we made an assumption about the parameter
 $q$ involved, which implies a restriction on $\beta = q+q^{-1}$.
In addition, we  assumed 
$A$ and $A^*$ are not nilpotent. With these comments in mind, 
we pose the following problem.

\begin{problem} Classify all the generalized 
TD pairs that are not TD pairs. 
\end{problem}

\noindent 
Given a field  $\fld$, and given scalars 
 $\beta, \gamma, \gamma^*,
\varrho, \varrho^*$ taken from $\fld$,
let 
 $T=T(\beta, \gamma, \gamma^*,
\varrho, \varrho^*)$  denote the associative $\fld$-algebra
with identity generated by two symbols $A$, $A^*$
subject to
the relations
(\ref{eq:qdolangradyfin1}), 
(\ref{eq:qdolangradyfin2}).
Let us call $T$ a {\it TD algebra}.
The generalized TD pairs are essentially the 
same thing as 
 irreducible finite dimensional
modules for  TD algebras. 
Given the results of this paper, we would  not
be suprised if the TD algebras had some  
interesting properties, and we encourage the
reader to look into this.

\medskip
\noindent 
We mention the relations 
(\ref{eq:qdolangradyfin1}),
(\ref{eq:qdolangradyfin2}) 
previously appeared in \cite{TersubIII}, in the 
context of 
$P$- and $Q$-polynomial association schemes
\cite{BanIto}, \cite{bcn},
\cite{Leopandq},
\cite{Tercharpq},
\cite{Ternew}.
In
 \cite{TersubIII}, it is shown 
 that for these schemes the 
matrices $A$ and $A^*$ in 
Example
\ref{ex:pandq}   
satisfy
(\ref{eq:qdolangradyfin1}), 
(\ref{eq:qdolangradyfin2}).
In this context the algebra 
generated by $A$ and $A^*$ is known as the subconstituent
algebra or the Terwilliger algebra
\cite{Cau},
\cite{Col},
\cite{Curbip12},
\cite{CurNom},
\cite{Go},
\cite{HobIto},
\cite{Tan},
\cite{TersubI},
\cite{TersubII}.

\medskip
\noindent 
We  mention the relations
(\ref{eq:qdolangradyfin1}),
(\ref{eq:qdolangradyfin2})  are satisfied by the
generators of 
both the classical
and quantum ``Quadratic  Askey-Wilson algebra''
introduced by   Granovskii, Lutzenko, and 
Zhedanov \cite{GYLZmut}.
See 
\cite{GYZTwisted},
\cite{GYZL}, 
\cite{SVVZPer},
\cite{ZheHid}, 
\cite{ZheHiggs},
\cite{ZheCart}
for more information on this algebra.

\medskip
\noindent Here are a few conjectures about TD systems.

\begin{conjecture}
\label{conj:boundonsubdeg}
Let $\Phi$ denote a TD system with diameter $d$, and let
$\rho_0, \rho_1, \ldots, \rho_d$ denote the associated scalars
from 
Corollary \ref{thm:symunimodal}.
We conjecture
\begin{equation}
\rho_i \;\leq \;\biggl( {{d}\atop {i}} \biggr) \qquad \qquad (0 \leq i \leq d).
\label{eq:upperbnd}
\end{equation}
\end{conjecture}
\noindent  The above conjecture has been proven for $d\leq 3$ by
Tanabe.

\begin{conjecture}
\label{con:Vspanguess}
Let $\fld $ denote a field, and 
let $V$ denote a vector space over $\fld$ with finite
positive dimension. Let 
$\Phi$  denote a TD system on $V$, and let $d$ denote
the diameter.
Let the maps 
$R$ and  $L$ be as in 
Definition
\ref{def:defRandL}.
Let the space $U_0$  be as in 
Definition
\ref{def:meaningofVi}, and 
pick any nonzero vector $v \in U_0$. We conjecture $V$ is spanned by
the vectors of the form
\beast
L^{i_1}R^{i_2}L^{i_3}R^{i_4} \cdots R^{i_n}v,
\eeast
where 
$i_1, i_2, \ldots, i_n$ ranges over all sequences
such that $n$ is nonnegative and even, and 
such that $i_1, i_2, \ldots i_n$ are integers
satisfying $0\leq i_1<i_2<\cdots <i_n\leq d$.
For example, for $d=3$, we are asserting
$V$ is spanned by
\beast
v,\quad  Rv, \quad R^2v, \quad  R^3v, 
\quad LR^2v, \quad LR^3v,
\quad L^2R^3v, \quad  
RL^2R^3v.
\eeast
\end{conjecture}
\noindent Conjecture 
\ref{con:Vspanguess}
has been proved for $d\leq 3$ by Tanabe. 
We remark the Conjecture 
\ref{con:Vspanguess}
implies 
Conjecture
\ref{conj:boundonsubdeg}.

\begin{conjecture}
\label{conjfactor}
Let $\Phi$ denote a TD system with diameter $d$, and let
$\rho_0, \rho_1, \ldots, \rho_d$ denote the corresponding
scalars from Corollary
\ref{thm:symunimodal}. Consider the polynomial
\begin{equation}
\sum_{i=0}^d \rho_i t^i
\label{eq:rhopoly}
\end{equation}
in a variable $t$. We conjecture
(\ref{eq:rhopoly}) equals  
\begin{equation}
(1+t+t^2+\cdots +t^{d_1})
(1+t+t^2+\cdots +t^{d_2})
\cdots
(1+t+t^2+\cdots +t^{d_n}),
\label{eq:seriesfactor}
\end{equation}
where 
$d_1, d_2, \ldots, d_n$ are positive integers
whose sum is $d$.
\end{conjecture}

\noindent We remark Conjecture 
\ref{conjfactor} implies
Conjecture \ref{conj:boundonsubdeg}.

\medskip
\noindent
Tatsuro Ito, Department of Computational Science, Faculty of Science,
Kanazawa University, Kakuma-machi, Kanazawa 920--1192, Japan \hfil\break
E-mail: ito@kappa.s.kanazawa-u.ac.jp \hfil\break

\medskip
\noindent 
Kenichiro Tanabe, Graduate School  of Mathematics, Kyushu University, 
33 Fukuoka 812-8581,
Japan \hfil \break
E-mail:  tanabe@math.kyushu-u.ac.jp \hfil\break

\medskip
\noindent
Paul Terwilliger, Department of Mathematics,
University of Wisconsin, 480 Lincoln drive, 
Madison, Wisconsin, 53706, USA \hfil\break
E-mail: terwilli@math.wisc.edu \hfil\break


\begin{thebibliography}{10}

\bibitem{CKOns}
C.~Ahn and K.~Shigemoto.
\newblock Onsager algebra and integrable lattice models.
\newblock {\em Modern Phys. Lett. A}, 6(38):3509--3515, 1991.

\bibitem{BanIto}
E.~Bannai and T.~Ito.
\newblock {\em Algebraic Combinatorics I: Association Schemes}.
\newblock Benjamin/Cummings, London, 1984.

\bibitem{bcn}
A.~E. Brouwer, A.~M. Cohen, and A.~Neumaier.
\newblock {\em Distance-Regular Graphs}.
\newblock Springer-Verlag, Berlin, 1989.

\bibitem{Cau}
J.~S. Caughman~{I}{V}.
\newblock The {T}erwilliger algebras of bipartite ${P}$- and ${Q}$-polynomial
  schemes.
\newblock {\em Discrete Math.}, 196(1-3):65--95, 1999.

\bibitem{CPqaa}
Vyjayanthi Chari and Andrew Pressley.
\newblock Quantum affine algebras.
\newblock {\em Comm. Math. Phys.}, 142(2):261--283, 1991.

\bibitem{CPqaar}
Vyjayanthi Chari and Andrew Pressley.
\newblock Quantum affine algebras and their representations.
\newblock In {\em Representations of groups (Banff, AB, 1994)}, pages 59--78.
  Amer. Math. Soc., Providence, RI, 1995.

\bibitem{Col}
B.~V.~C. Collins.
\newblock The girth of a thin distance-regular graph.
\newblock {\em Graphs Combin.}, 13(1):21--30, 1997.

\bibitem{Curbip12}
B.~Curtin.
\newblock Bipartite distance-regular graphs, parts {I} and {II}.
\newblock {\em Graphs Combin.}, to appear.

\bibitem{CurNom}
B.~Curtin and K.~Nomura.
\newblock Distance-regular graphs related to the quantum enveloping algebra of
  $sl(2)$.
\newblock {\em J. Algebraic Combin.}, to appear.

\bibitem{DateRoan}
E.~Date and S.S. Roan.
\newblock The algebraic structure of the onsager algebra.
\newblock {\em Czech. J. Physics}.

\bibitem{Dav}
B.~Davies.
\newblock Onsager's algebra and the {D}olan-{G}rady condition in the
  non-self-dual case.
\newblock {\em J. Math. Phys.}, 32(11):2945--2950, 1991.

\bibitem{Dolgra}
L.~Dolan and M.~Grady.
\newblock Conserved charges from self-duality.
\newblock {\em Phys. Rev. D (3)}, 25(6):1587--1604, 1982.

\bibitem{Go}
J.~Go.
\newblock The {T}erwilliger algebra of the {H}ypercube ${Q_D}$.
\newblock {\em Graphs Combin.}, submitted.

\bibitem{GYLZmut}
Ya.~I. Granovski{\u\i}, I.~M. Lutzenko, and A.~S. Zhedanov.
\newblock Mutual integrability, quadratic algebras, and dynamical symmetry.
\newblock {\em Ann. Physics}, 217(1):1--20, 1992.

\bibitem{GYZTwisted}
Ya.~I. Granovski{\u\i} and A.~S. Zhedanov.
\newblock ``{T}wisted'' {C}lebsch-{G}ordan coefficients for ${\rm {s}{u}}\sb
  q(2)$.
\newblock {\em J. Phys. A}, 25(17):L1029--L1032, 1992.

\bibitem{GYZL}
Ya.~I. Granovski{\u\i}, A.~S. Zhedanov, and I.~M. Lutsenko.
\newblock Quadratic algebras and dynamical symmetry of the {S}chr\"odinger
  equation.
\newblock {\em Soviet Phys. JETP}, 72(2):205--209, 1991.

\bibitem{HobIto}
S.~Hobart and T.~Ito.
\newblock The structure of nonthin irreducible ${T}$-modules of endpoint 1:
  ladder bases and classical parameters.
\newblock {\em J. Algebraic Combin.}, 7(1):53--75, 1998.

\bibitem{Leopandq}
D.~A. Leonard.
\newblock Parameters of association schemes that are both ${P}$- and
  ${Q}$-polynomial.
\newblock {\em J. Combin. Theory Ser. A}, 36(3):355--363, 1984.

\bibitem{Per}
J.~H.~H. Perk.
\newblock Star-triangle equations, quantum {L}ax pairs, and higher genus
  curves.
\newblock In {\em Theta functions---Bowdoin 1987, Part 1 (Brunswick, ME,
  1987)}, pages 341--354. Amer. Math. Soc., Providence, RI, 1989.

\bibitem{Roanmpi}
S.S. Roan.
\newblock Onsager's algebra, loop algebras and chiral potts model.
\newblock 70, 1991.

\bibitem{SVVZPer}
V.~Spiridonov, L.~Vinet, and A.~Zhedanov.
\newblock Periodic reduction of the factorization chain and the {H}ahn
  polynomials.
\newblock {\em J. Phys. A}, 27(18):L669--L675, 1994.

\bibitem{Tan}
K.~Tanabe.
\newblock The irreducible modules of the {T}erwilliger algebras of {D}oob
  schemes.
\newblock {\em J. Algebraic Combin.}, 6(2):173--195, 1997.

\bibitem{Tercharpq}
P.~Terwilliger.
\newblock A characterization of ${P}$- and ${Q}$-polynomial association
  schemes.
\newblock {\em J. Combin. Theory Ser. A}, 45(1):8--26, 1987.

\bibitem{TersubI}
P.~Terwilliger.
\newblock The subconstituent algebra of an association scheme. {I}.
\newblock {\em J. Algebraic Combin.}, 1(4):363--388, 1992.

\bibitem{TersubII}
P.~Terwilliger.
\newblock The subconstituent algebra of an association scheme. {I}{I}.
\newblock {\em J. Algebraic Combin.}, 2(1):73--103, 1993.

\bibitem{TersubIII}
P.~Terwilliger.
\newblock The subconstituent algebra of an association scheme. {I}{I}{I}.
\newblock {\em J. Algebraic Combin.}, 2(2):177--210, 1993.

\bibitem{Ternew}
P.~Terwilliger.
\newblock A new inequality for distance-regular graphs.
\newblock {\em Discrete Math.}, 137(1-3):319--332, 1995.

\bibitem{LS99}
P.~Terwilliger.
\newblock Two linear transformations each tridiagonal with respect to an
  eigenbasis of the other.
\newblock {\em Linear Algebra Appl}, submitted.

\bibitem{Ugl}
D.~B. Uglov and I.~T. Ivanov.
\newblock ${\rm sl}({N})$ {O}nsager's algebra and integrability.
\newblock {\em J. Statist. Phys.}, 82(1-2):87--113, 1996.

\bibitem{ZheHid}
A.~S. Zhedanov.
\newblock ``{H}idden symmetry'' of {A}skey-{W}ilson polynomials.
\newblock {\em Teoret. Mat. Fiz.}, 89(2):190--204, 1991.

\bibitem{ZheHiggs}
A.~S. Zhedanov.
\newblock The ``{H}iggs algebra'' as a ``quantum'' deformation of ${\rm
  {s}{u}}(2)$.
\newblock {\em Modern Phys. Lett. A}, 7(6):507--512, 1992.

\bibitem{ZheCart}
A.~S. Zhedanov.
\newblock Quantum ${\rm {s}{u}}\sb q(2)$ algebra: ``{C}artesian'' version and
  overlaps.
\newblock {\em Modern Phys. Lett. A}, 7(18):1589--1593, 1992.

\end{thebibliography}
\end{document}